\documentclass[12pt]{article}
\usepackage{latexsym,amssymb,color}
\usepackage{epsfig}
\newcommand{\eq}{\begin{equation}}
\newcommand{\en}{\end{equation}}
\newcommand{\re}[1]{\mbox{(\ref{#1})}}

\newtheorem{Theorem}{Theorem}
\newtheorem{theorem}[Theorem]{Theorem}
\newtheorem{lemma}[Theorem]{Lemma}
\newtheorem{corollary}[Theorem]{Corollary}
\newtheorem{construction}[Theorem]{Construction}
\newtheorem{proposition}[Theorem]{Proposition}
\newtheorem{example}[Theorem]{Example}
\newtheorem{defn}[Theorem]{Definition}

\newtheorem{question}[Theorem]{Question}
\newtheorem{conjecture}[Theorem]{Conjecture}
\newtheorem{condition}[Theorem]{Condition}
\newtheorem{remark}[Theorem]{Remark}
\newtheorem{problem}[Theorem]{Problem}
\newtheorem{jpfigure}[Theorem]{Figure}

\def\proof{\noindent{\bf Proof.\ \ }}

\def\endpf{\hfill $\Box$ \vskip .25in}

\newfont{\msbm}{msbm10 at 12pt}
\newfont{\eusb}{eusb10}
\newfont{\eusm}{eusm10}
\newfont{\eurb}{eurb10}
\newfont{\eurm}{eurm10}
\newfont{\eufb}{eufb10}
\newfont{\eufm}{eufm10}

\newcommand {\ints} {\mbox{\msbm\symbol{'132}}}

\newcommand{\con}{\rightarrow}

\newcommand{\ed}{\mbox{$ \ \stackrel{d}{=}$ }}

\newcommand{\eps}{\varepsilon}

\newcommand{\Tak}{T{\'a}kacs}

\newcommand{\bbc}[2]{\left( \hspace{-.5em}{#1\choose
      #2}\hspace{-.5em}\right)} 
\newcommand{\tbbc}[2]{\left(\hspace{-.3em}\left({#1\atop
        #2}\right)\hspace{-.3em}\right)} 
\newcommand{\be}{\begin{enumerate}}
\newcommand{\ee}{\end{enumerate}}
\newcommand{\beq}{\begin{equation}}
\newcommand{\eeq}{\end{equation}}
\newcommand{\bea}{\begin{eqnarray}}
\newcommand{\eea}{\end{eqnarray}}
\newcommand{\beas}{\begin{eqnarray*}}
\newcommand{\eeas}{\end{eqnarray*}}
\newcommand{\bm}[1]{{\mbox{\boldmath $#1$}}}
\newcommand{\sbm}[1]{{{{\mbox{${\mbox{\scriptsize\boldmath$#1$}}$}}}}}
\newcommand{\suu}{{{{\mbox{${\mbox{\scriptsize\boldmath$u$}}$}}}}}
\newcommand{\tuu}{{{{\mbox{${\mbox{\tiny\boldmath$u$}}$}}}}}
\newcommand{\pa}{\mathrm{park}} 
\newcommand{\nreals}{\mathbb{R}_{\geq 0}}
\newcommand{\reals}{\mathbb{R}}

\newcommand{\lb}[1]{\label{#1}}

\newenvironment{thm}[1]{\begin{theorem}\label{#1}}{\end{theorem}}

\newenvironment{crl}[1]{\begin{corollary}\protect\label{#1}}{\end{corollary}}

\newenvironment{exm}[1]{\begin{example}\protect\label{#1}}{\end{example}}

\newcommand{\xx}{{\mbox{\boldmath$x$}}}
\newcommand{\yy}{{\mbox{\boldmath$y$}}}
\newcommand{\zz}{{\mbox{\boldmath$z$}}}
\newcommand{\uu}{{\mbox{\boldmath$u$}}}
\newcommand{\rr}{{\mbox{\boldmath$r$}}}
\newcommand{\vzero}{{\mbox{\boldmath$0$}}}
\newcommand{\vone}{{\mbox{\boldmath$1$}}}
\renewcommand{\ss}{{\mbox{\boldmath$s$}}}
\newcommand{\kk}{{\mbox{\boldmath$k$}}}

\newcommand{\skk}{{{{\mbox{${\mbox{\scriptsize\boldmath$k$}}$}}}}}

\newcommand{\hN}{\widehat{N}}

\begin{document}
\title{\textcolor{red}{A polytope related to empirical distributions,
    plane trees,  
parking functions, and the associahedron }}
\author{by Jim Pitman\thanks{Research supported in part by NSF grant
    97-03961}\ \ and Richard Stanley\thanks{Research supported in part by
    NSF grant 95-00714}\\ 
\\
Technical Report No. 560
\\
\\
Department of Statistics\\
University of California\\
367 Evans Hall \# 3860\\
Berkeley, CA 94720-3860\\[.1in]
{\small version of 8 July 1999}}

\maketitle

\renewcommand{\insert}{{\rm insert}}
\newcommand{\tent}{{\rm tent}}
\newcommand{\giv}{{\,|\,}}
\newcommand{\tree}{{\rm TREE}}
\newcommand{\TREE}{{\rm TREE}}
\newcommand{\hC}{\widehat{C}}
\renewcommand{\eps}{\varepsilon}
\newcommand{\epstd}{\varepsilon^*}
\newcommand{\teps}{\widetilde{\eps}}
\newcommand{\sU}{U^\dagger}
\newcommand{\skU}{U}
\newcommand{\preals}{\reals_+}
\newcommand{\nnints}{\mathbb{N}}
\newcommand{\sV}{V}
\newcommand{\tC}{\widetilde{C}}

\newcommand{\thv}{{\mbox{\boldmath$\theta$}}}
\newcommand{\sthv}{{{{\mbox{${\mbox{\scriptsize\boldmath$\theta$}}$}}}}}
\newcommand{\bcth}{{\mbox{\boldmath$\theta$}}}
\newcommand{\thbc}{{\mbox{\boldmath$\theta$}}}
\newcommand{\tbc}{{ \theta \bc}}
\newcommand{\bl}{{\mbox{\boldmath$l$}}}
\newcommand{\el}{l}
\newcommand{\bt} { {\bf t }}
\newcommand{\dd} {n}
\newcommand{\BT} { {\bf T }}
\newcommand{\BTn} { {\bf B }_n}
\newcommand{\BTT} { \widetilde{{\bf B }}}
\newcommand{\TT} { {{\cal T }}}
\newcommand{\TTT} { \widetilde{{\cal T }}}
\newcommand{\BB} { \widetilde{{\cal B }}}
\newcommand{\BBd} {{ \cal  B}}
\newcommand{\YY} {{ \cal Y}}
\newcommand{\GG} {{ \cal G}}
\newcommand{\SV} { \SS_\bc^V }
\newcommand{\hubs} { \mbox{\rm hubs}}
\newcommand{\Xo} { X_1 }
\newcommand{\Xt} { X_2 }
\newcommand{\Xth} { X_3 }
\newcommand{\Xn} { X_n }
\newcommand{\tJ} { J_3 }
\newcommand{\Hn} { { \cal H } _ n}
\newcommand{\Hth} { { \cal H } _ 3}
\newcommand{\pitil} {{ \tilde{{\pi}}  }}
\newcommand{\sfrac}[2]{{\textstyle\frac{#1}{#2}}}
\newcommand{\tail} { {\rm tail}}
\newcommand{\tails} { {\rm tails}}
\newcommand{\btt} { \widetilde{{\bf t }}}
\newcommand{\btth} { \widetilde{{\bf t }}^*}
\newcommand{\PP} { {\cal P }}
\newcommand{\Lp} { {L_{1^+,2^+}}}
\newcommand{\PPT} { \tilde{{\cal P }}}
\newcommand{\lic} { { L_{I,c} } }
\newcommand{\fic} { { f_{I,c} (s)} }
\newcommand{\bnk} {B_{n,k}}
\newcommand{\xell} {x}
\newcommand{\la} { \lambda}
\newcommand{\hf}{\sfrac{1}{2}}
\newcommand{\sxth}{\sfrac{1}{6}}
\newcommand{\ssc} { {\cal S _{{\bf c} }}}
\newcommand{\TTc} {  {\cal T _{{\bf c} }}}
\newcommand{\TTH} {  \widehat{{\cal T}}}
\newcommand{\ssch} [1]{ { {{\cal S}} ^{(#1)} _{{\bf c} }} }
\newcommand{\sscn} { {\cal S _{{ {\bf c} {\it (n)}} }}}
\newcommand{\ssco} { {\cal S _{{ {\bf c} {\it (1)} } }}}
\newcommand{\met} [1] { {\rm Tree}_{#1} }
\renewcommand{\span} [1] { {\rm Span}_{#1} }
\newcommand{\SI} { S_I }
\newcommand{\SIu} { S_{I,u} }
\newcommand{\SID} { S_I ^\Delta}
\newcommand{\SnpD} { S_{n+1} ^\Delta}
\renewcommand{\SS} { {\cal S }}
\newcommand{\schr} { {\cal S }}
\newcommand{\cc} { {\bf c }}
\newcommand{\shape} {\mbox{\rm shape}}
\newcommand{\ned} { \# {\cal E }}
\newcommand{\exd} { D}
\newcommand{\flr}[1]{ { \lfloor #1 \rfloor } }
\newcommand{\bc}{{\bf c}}
\newcommand{\bu}{{\bf u}}
\newcommand{\br}{{\bf r}}
\newcommand{\bs}{{\bf s}}
\newcommand{\bU}{{\bf U}}
\newcommand{\RR}{{\cal R}}
\newcommand{\UU}{{\cal U}}
\newcommand{\CC}{{\cal C}}
\newcommand{\ecs}{e^{-\frac{1}{2}s^2 - sc}}

\newcommand{\bP}{{\bf P}}
\newcommand{\bT}{{\bf T}}
\newcommand{\bZ}{{\bf Z}}
\newcommand{\bJ}{{\bf J}}
\newcommand{\bQ}{{\bf Q}}
\newcommand{\bX}{{\bf X}}
\newcommand{\bY}{{\bf Y}}
\newcommand{\bV}{{\bf V}}
\newcommand{\bC}{{\bf C}}
\newcommand{\ED}{\stackrel{\con}{{\cal E}} }
\newcommand{\var}{{\rm var}\ }
\newcommand{\Ag}{||g||_2^2}
\begin{abstract}
The volume of the $\dd$-dimensional 
polytope
$$\Pi_\dd(\xx):= \{ \yy \in \reals^\dd: y_i \ge 0 \mbox{ and }  y_1 +
\cdots+ y_i \le x_1 + \cdots + x_i \mbox{ for all } 1 \le i \le \dd \}
$$
for arbitrary $\xx:=(x_1, \ldots, x_\dd)$ with $x_i >0$ for all $i$ defines
a polynomial in variables $x_i$ which admits a number of interpretations,
in terms of empirical distributions, plane partitions, and
parking functions. We interpret the terms of this polynomial 
as the volumes of chambers in two different polytopal subdivisions of
$\Pi_\dd(\xx)$. The first of these subdivisions generalizes to
a class of polytopes called sections of order cones.
In the second subdivision, the chambers are indexed in a natural way by
rooted binary trees with $n+1$ vertices, and the configuration of these
chambers provides a representation of another polytope with many applications, 
the {\em associahedron}.
\end{abstract}



\section{Introduction}

\medskip
The focal point of this paper is the $\dd$-dimensional 
polytope
$$\Pi_\dd(\xx):= \{ \yy \in \reals^\dd: y_i \ge 0 \mbox{ and }  y_1 +
\cdots+ y_i \le x_1 + \cdots + x_i \mbox{ for all } 1 \le i \le \dd \}
$$
for arbitrary $\xx:=(x_1, \ldots, x_\dd)$ with $x_i >0$ for all $i$.
The $\dd$-dimensional volume 
$$
V_\dd(\xx):= \mathrm{Vol} (\Pi_\dd(\xx))$$
is a homogeneous polynomial of degree $\dd$ in the variables
$x_1, \ldots, x_\dd$, which we call the \emph{volume polynomial}.
This polynomial arises naturally in several different settings:
in the calculation of probabilities
derived from empirical distribution functions or the order statistics
of $\dd$ independent random variables (see \S \ref{edfs}),
and in the study of parking functions and plane partitions 
(see \S \ref{sec:pfpp}).
See also Marckert and Chassaing \cite{marchas99} regarding similar
connections between the theories of parking functions, empirical processes, 
and rooted trees.

\begin{figure}
\centerline{\epsfig{figure=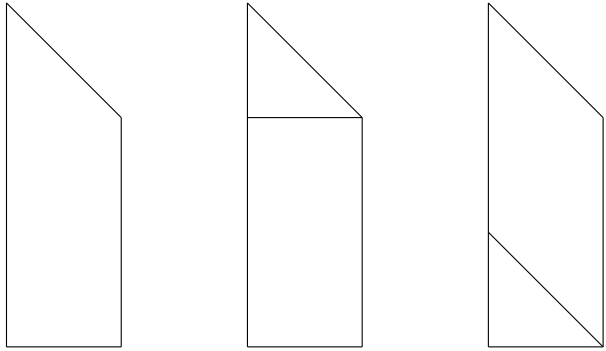,height=2in}}
\caption{$\Pi_2(\xx)$ and its two subdivisions}
\label{fig2}
\end{figure}

Trivially, $V_1(\xx) = x_1$.
The formula 
$$V_2(\xx) =x_1 x_2 + \hf x_1^2$$
has two
natural interpretations by a subdivision of
$\Pi_2(\xx)$ into 2 pieces of areas $x_1 x_2$ and $\hf x_1^2$, as shown
in Figure~\ref{fig2} 
for horizontal coordinate $x_1 =1$ and vertical coordinate $x_2 =2$.

\begin{figure}
\centerline{\epsfig{figure=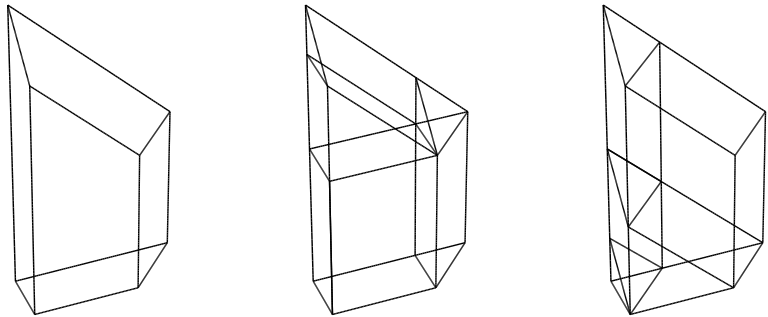,height=2in}}
\caption{$\Pi_3(\xx)$ and its two subdivisions}
\label{fig3}
\end{figure}

The 5 terms of 
\eq
\lb{v3}
V_3(\xx) = x_1 x_2 x_3 + \hf  x_1^2 x_2 + \hf x_1 x_2^2 + \hf x_1^2 x_3 + \sxth x_1^3
\en
can be interpreted in two ways as the volumes determined by 
two different subdivisions of $\Pi_3(\xx)$ into 5 chambers,
as in the perspective diagrams of Figure~\ref{fig3}
where $x_i = i$ for $i = 1,2,3$, 
the first coordinate points out of the page, the second to the right
and the third up, and the viewpoint is $(5,-2,4)$.

A central result of this paper is the general formula 
for the volume polynomial which we present in the following theorem. 
Section \ref{edfs} offers a simple probabilistic proof of this theorem.
We show in Section \ref{pinsec} 
how this argument can also be interpreted geometically by a subdivision of 
$\Pi_\dd(\xx)$ into a collection of $\dd$-dimensional chambers, with
the volume of each chamber corresponding to a term of the volume
polynomial. This generalizes the subdivisions of $\Pi_2$ and $\Pi_3$ 
shown in the right hand panels of Figures ~\ref{fig2} and~\ref{fig3}.
Technically, by a \emph{subdivision} of $\Pi_\dd(\xx)$ we mean a
\emph{polytopal subdivision} in the sense of Ziegler \cite[p.\
129]{ziegler95}, and we call the $\dd$-dimensional polytopes involved
the \emph{chambers} of the subdivision.
The subdivision of $\Pi_\dd(\xx)$ described in
Section \ref{pinsec} is a specialization of a result 
presented in Section~\ref{ordercone} in the general context of 
``sections of order cones''. 
Section \ref{sec:assoc} shows how the subdivisions shown 
in the left hand panels of Figures ~\ref{fig2} and~\ref{fig3}
can be generalized to arbitrary $\dd$.  The chambers of this subdivision
of $\Pi_\dd(\xx)$ are indexed in a natural way by rooted binary plane
trees with $\dd +1$ leaf vertices, and the configuration of these
chambers provides a representation of another interesting polytope
with many applications, known as the \emph{associahedron}. 

\begin{thm} {volume}
For each $\dd = 1,2, \ldots$,
\eq
\lb{vform}
V_\dd(\xx) = \sum_{\skk \in K_\dd} \,
\prod_{i=1}^\dd {x_i ^ {k_i} \over k_i!} = \frac{1}{\dd!}
  \sum_{\skk \in K_\dd} {\dd\choose k_1,\dots,k_\dd}x_1^{k_1}\cdots
   x_\dd^{k_\dd},
\en
where 
\eq
\lb{kdef}
K_\dd:= \{ \kk \in \nnints^\dd: \sum_{i=1}^j k_i \ge j \mbox{ \rm for all }
1 \le i \le \dd-1  \mbox{ \em and }
\sum_{i=1}^\dd k_i = \dd \}
\en
with $\nnints := \{0,1,2, \ldots \}$.
\end{thm}

In particular, the number of nonzero coefficients in $V_\dd$ is the
number of elements of $K_\dd$, which is well known to be the $\dd$th
Catalan number $C_\dd$ (see e.g.\ \cite[Exer.\ 6.19(w)]{stanleyv2} for
a simple variant), the first few of which are $1,2,5,14,42,132,
\ldots$:
\eq
\lb{numt}
\# K_\dd = C_\dd := {1 \over \dd+1} { 2 \dd \choose \dd}.
\en

Formula \re{vform} should be compared with the following
alternate formula, which as indicated in Section \ref{edfs} can
be read from a formula of Steck \cite{steck69,steck71} for 
the cumulative distribution function of the random vector of order statistics
of $\dd$ independent random variables with uniform distribution on an
interval:
\eq
\lb{steck}
V_\dd(\xx) = \det \left[ 
{1 (j-i+1 \ge  0) \over (j-i+1)! } \left( \sum_{h=1}^i x_h \right)^{j-i+1} 
\right]_{1 \le i,j \le n}
\en
where $\det \left[ a_{ij} \right] _{1 \le i,j \le n}$ denotes the determinant
of the $\dd \times n$ matrix with entries $a_{ij}$,
and $1(\cdots)$ equals $1$ if $\cdots$ and $0$ else.
See \cite{ejgp72} for an elementary probabilistic proof of \re{steck}.
This formula allows the
expansion of $V_\dd(\xx)$ into monomial terms to be generated for
arbitary $\dd$ by just few lines of \emph{Mathematica} code. 

Another formula of Steck \cite{steck69,steck71}, with an elementary proof 
in \cite{ejgp72}, gives 
the number $\#(b,c)$ of $j \in \ints^\dd$ with 
$j_1 < j_2 < \cdots < j_\dd$ and $b_i < j_i < c_i$
for all $1 \le i \le n$ for arbitrary $b,c \in \ints^\dd$ with
$b_1 \le b_2 \le  \cdots < b_\dd$ and $c_1 \le c_2 \le  \cdots < c_\dd$:
\eq
\lb{steck2}
\# (b,c) = \det \left[ 
{1 (j-i+1 \ge  0, c_i - b_j >1)  } 
{ c_i -b_j + j - i - 1 \choose j-i +1 }
\right]_{1 \le i,j \le n} .
\en
We explain after the proof of Theorem~\ref{pp} how
these formulae (\ref{steck}) and (\ref{steck2}) can be deduced
from a result of MacMahon on the enumeration of plane partitions.

In Section \ref{edfs} we deduce the following special evaluations
of the volume polynomial from some well known results in the theory
of empirical distributions: for $a,b \ge 0$
\eq
\lb{abform}
n! V_\dd( a , b, \ldots, b ) = a ( a  + n b )^{\dd-1}
\en
while for $\dd \ge 3$ and $a,b,c \ge 0$
\eq
\lb{abc}
n! V_\dd( a , \overbrace{ b, \ldots, b}^{n - 2 \rm \,places },c) = 
a ( a  + n b )^{n-1} + n a (c - b ) ( a  + (n-1) b )^{n-2}
\en
and for $\dd \ge 3$, $1 \le m \le n -2$ and $a,b,c \ge 0$
\eq
\lb{abcd1}
n! V_\dd( a , \overbrace{ b, \ldots, b,}^{n - m - 1 \rm \, places }c, \,
\overbrace{ 0, \ldots, 0}^{m - 1 \rm \, places })
= a \, \sum_{ j = 0}^ m { n \choose j } ( c - ( m+1 - j ) b )^j (a + (n-j) b )^{n-j-1} .
\en
As we indicate in Section \ref{sec:pfpp}, these formulae read from the
theory of empirical distributions have interesting
combinatorial interpretations in terms of parking functions and 
plane partitions.

\section{Uniform Order Statistics and Empirical Distribution Functions}
\label{edfs}

Let $(U_{\dd,i}, 1 \le i \le n)$ be the \emph{order statistics} of $\dd$
independent uniform $(0,1)$ variables $U_1, U_2, \ldots, U_\dd$.  That
is to say, $U_{\dd,1} \le U_{\dd,2} \le \cdots \le U_{\dd,n}$ are the ranked
values of the $U_i, 1\le i \le n$.  
Because the random vectors 
$(U_{\dd,j}, 1 \le j \le n)$ and
$(1 - U_{\dd,n+1-j}, 1 \le j \le n)$ 
have the same uniform distribution with constant density $\dd!$ on the simplex
\eq
\lb{simplex}
\{ \uu \in \reals^\dd : 0 \le u_1 \le \cdots \le u_\dd \le 1 \}
\en
for arbitrary vectors $\rr$ and $\ss$ in this simplex there are
the formulae
\eq
\lb{os}
P( U_{\dd,j}  \le s_j \mbox{ for all } 1 \le j \le n) =  
n! \, V_\dd(x_1, \ldots, x_\dd ) \mbox{ where } x_j:= s_j - s_{j-1}
\en
where $s_0:= 0$ and
\eq
\lb{os1}
P( U_{\dd,j}  \ge  r_j \mbox{ for all } 1 \le j \le n) 
=  n! \,
V_\dd(x_1, \ldots, x_\dd ) \mbox{ where } x_j:= r_{\dd+2-j} - r_{\dd+1-j}
\en
where $r_{\dd+1}:= 1$.
Thus the probability 
\eq
\lb{os2}
P_\dd(\rr,\ss):= P( r_j \le U_{\dd,j}  \le s_j \mbox{ for all } 1 \le
j \le n)  
\en
can be evaluated in terms of $V_\dd$ if either $\rr = \vzero$ or $\ss = \vone$.
See \cite[\S 9.3]{sw86} for a review of results involving
these probabilities, including various recursion formulae which are
useful for their computation.

\noindent
{\bf Proof of Theorem \ref{volume}}.
By homogeneity of $V_{\dd}$,
it suffices to prove the formula when $s_\dd \le 1$.
Fix $\xx$ and consider the probability \re{os}.
For $1 \le i \le n+1$ let $N_i$ denote the number of $U_{\dd,i}$ that
fall in the interval $(s_{i-1},s_i]$, with the conventions $s_0 = 0$
and $s_{\dd+1} =1$:  
\eq
\lb{ndef}
N_i := \sum_{i=1}^\dd 1( s_{i-1} < U_{\dd,j}  \le s_i) = 
\sum_{i=1}^\dd 1( s_{i-1} < U_j  \le s_i).
\en
The second expression for $N_i$ shows that
the random vector $(N_i, 1 \le i \le n+1)$ has the
\emph{multinomial distribution with parameters $\dd$ 
and } $(x_1, \ldots, x_\dd, x_{\dd+1})$ for $x_{i}:= s_{i}- s_{i-1}$,
meaning that for each vector of $\dd+1$ nonnegative integers 
$(k_i, 1 \le i \le n+1)$ with
$\sum_{i=1}^{\dd+1} k_i = n$, we have
\eq
\lb{multi1}
P( N_i = k_i, 1 \le i \le n+1) = 
n! \prod_{i=1}^{\dd+1}  { {x_i}^{k_i} \over k_i !} .
\en
By definition of the $U_{\dd, j}$ and \re{ndef}, 
the events $(U_{n,j} \le s_j)$ and 
$( \Sigma_{i=1}^j N_i \ge j )$ are identical.
Thus
$$
P( U_{n,j}  \le s_j \mbox{ for all } 1 \le j \le n)
= P( \Sigma_{i=1}^j N_i \ge j \mbox{ for all } 1 \le j \le n)
$$
$$
= \sum_{\skk \in K_\dd} P( N_i = k_i, 1 \le i \le n, N_{n+1} = 0) 
= n! \, \sum_{\skk \in K_\dd} \prod_{i=1}^{n}  { {x_i}^{k_i} \over k_i !}
$$
by application of \re{multi1} with $k_{n+1} = 0$.
Compare the result of this calculation with \re{os} to obtain
\re{vform}.
\endpf

It is easily seen that the decomposition of 
the event \re{os} considered in the above argument corresponds to
a polytopal subdivision of $\Pi_\dd(\xx)$ which for $\dd=2$ and $\dd=3$ is that 
shown in the right hand panels of Figures ~\ref{fig2} and~\ref{fig3}.
See Section \ref{pinsec} for further discussion of this subdivision of 
$\Pi_\dd(\xx)$.

The following corollary of Theorem \ref{volume}
spells out two more probabilistic interpretations of $V_\dd$.

\begin{crl}{multi}
Let $(N_i, 1 \le i \le n+1)$ be a random vector with
multinomial distribution with parameters $\dd$ and 
$(p_1, \ldots, p_{n+1})$, as if
$N_i$ is the number of times $i$ appears 
in a sequence of $\dd$ independent trials with probability $p_i$
of getting $i$ on each trial for $1 \le i \le n+1$, where
$\sum_{i=1}^{n+1} p_i =1$.
Then
\eq
\lb{multiprob}
P( \Sigma_{j=1}^i N_j \ge i \mbox{ for all } 1 \le i \le n ) = 
n! \,V_\dd(p_1, p_2, \ldots, p_\dd) .
\en
and
\eq
\lb{multiprob2}
P( \Sigma_{j=1}^i N_j < i \mbox{ for all } 1 \le i \le n ) = 
n! \,V_\dd(p_{n+1}, p_\dd, \ldots, p_2) .
\en
\end{crl}
\proof
The first formula is read from the previous proof of \re{vform}.
The second is just the first applied to
$
(\hN_{1}, \ldots, \hN_{n+1}):= (N_{n+1}, \ldots, N_1)
$
instead of $(N_{1}, \ldots, N_{n+1})$, because
$$
\sum_{i=1}^j \hN_i = \sum_{i=1}^j N_{n+2 - i} = n - \sum_{i=1}^{n+1-j} N_i
$$
so that 
$$
\sum_{i=1}^j \hN_i \ge j \mbox{ iff } \sum_{i=1}^{n+1-j} N_i <  n +1 - j,
$$
and hence the event that $\sum_{i=1}^j \hN_i \le  j $ for all 
$1 \le j \le n$ is identical to the event that $\sum_{i=1}^{m}N_i < m$ for 
all $1 \le m \le n$.
\endpf

Let
$$
F_\dd(t) := {1 \over n } \sum_{i =1}^\dd 1 ( U_i \le t) = {1 \over n }
\sum_{i =1}^\dd 1 ( U_{\dd,i} \le t) 
$$
be the usual \emph{empirical distribution function} associated with
the uniform random sample $U_1, \ldots, U_\dd$. So $F_\dd$ rises by a step 
of $1/n$ at each of the sample points.
It is well known \cite{sw86} that for
any for continuous increasing functions $f$ and $g$, the probability
$$
P( f(t) \le F_\dd(t) \le g(t) \mbox{ for all } t )
$$
equals
$P_\dd(\rr,\ss)$ as in \re{os2}
where $\rr$ and $\ss$ are easily expressed in terms of
values of the inverse functions of $f$ and $g$ at $i/n$ for $0 \le i \le n$.
As an example, Daniels \cite{daniels45} discovered the remarkable
fact that for $0 \le p \le 1$ the probability
that the empirical distribution function does not cross the
line joining $(0,0)$ to $(p,1)$ equals $1-p$, no matter what 
$\dd = 1,2, \ldots$:
\eq
\lb{daniels}
P( F_\dd(t) \le t/p \mbox{ for all } 0 \le t \le 1 ) = 1 - p 
\en
which can be rewritten as
\eq
\lb{daniels1}
P( U_{\dd ,i} \ge i p/n \mbox{ for all } 1 \le i \le n) = 1 - p .
\en
As observed in \cite[Chapter X]{pitman79inf}, Daniels' formula
\re{daniels} can be understood 
without calculation by an argument which gives the stronger
result of \Tak\ \cite[Theorem 13.1]{takacs77}
that this formula holds with $F_\dd$ replaced by $F$ for
any random right-continuous non-decreasing step function $F$ with 
cyclically exchangeable increments and $F(0)=0$ and $F(1) =1$.
Essentially, this is a continuous parameter form of the ballot theorem.
Many other proofs of Daniels' formula are known: see 
\cite[\S 9.1]{sw86} and papers cited there. 
The form \re{daniels1} of Daniels' formula is equivalent via
\re{os1} to 
\eq
\lb{vnp}
n! V_\dd( 1-p, p/n, \ldots, p/n ) = 1 - p  
\en
for $0 \le p \le 1$.  By homogeneity of $V_\dd$, this amounts to the 
identity \re{abform} of polynomials in two variables $a$ and $b$.

Pyke \cite[Lemma 1]{pyke59} found the following 
formula: for all real $b$ and $x$ with 
\eq
\lb{bxcons}
0 \le b \le 1 \mbox{ and } 0 \le n b - x \le 1,
\en
\eq
P \left( \max_{1 \le i \le n} ( b i - U_{n,i} ) \le x \right)
=
\lb{pyke2}
(1 + x - n b ) \sum_{ j = 0 }^ {\lfloor x/a \rfloor} {n \choose j} ( j b - x )^j ( 1 + x - j b )^{n - j - 1} .
\en
As indicated in \cite[p. 354, Exercise 2]{sw86}, this formula gives
an expression for the probability that the empirical cumulative 
distribution function based 
on a sample of $\dd$ independent uniform $(0,1)$ variables crosses an
arbitrary straight line through the unit square.
See \cite[\S 9.1]{sw86} for proof of an equivalent of \re{pyke2},
various related results, and further references.
The identity in distribution
$$
(U_{n,i}, 1 \le i \le n ) \ed (1- U_{n,n+1-i},1 \le i \le n)
$$
shows that the probability in \re{pyke2} equals
\eq
\lb{pyke3}
P ( U_{n,i} \le 1 +x - n b + b ( i - 1) \mbox{ for all } 1 \le i \le n)
\en
which according to \re{os} is equal in turn to
\eq
\lb{pyke4}
n! V_\dd (x_1, \ldots , x_\dd ) \mbox{ for } 
x_i = \left\{
\begin{array}{lll}
1 + x - n b && \mbox{ if } i = 1 \\
b && \mbox{ if } 2 \le i < n - \lfloor x/a \rfloor + 1 \\
(n-i + 2 ) b - x && \mbox{ if } i =  n - \lfloor x/a \rfloor + 1 \\
0 && \mbox{ if } i >  n - \lfloor x/a \rfloor + 1.
\end{array}
\right.
\en
For $a := 1 +x - n b$ and $b$ subject to \re{bxcons}, that is
$0 < a \le 1$ and $0 \le b \le 1$,
the above discussion gives us equality of
\re{pyke2} and \re{pyke4} with $x = a + n b - 1$.
In particular,
provided $0 \le x < a$ there is only a term for $j=0$ in \re{pyke2},
so the equality of \re{pyke2} and \re{pyke4} reduces to \re{abform}.
Similarly, for $a \le x < 2a$ there are only terms for $j=0$ and $j=1$ in 
\re{pyke2}.
For $\dd \ge 3$ this allows us to deduce \re{abc} from \re{pyke2} 
first for $a,b,c >0$ with $a + (n-2) b + c =1$ and $ c < b$, 
thence as an identity of polynomials in $a,b,c$.
Similarly, for $\dd \ge 3$ and $1 \le m \le n -2$
when $\lfloor{x/a} \rfloor = m$ we obtain
the identity \re{abcd1} of polynomials in $a,b,c$.

According to Steck \cite{steck69,steck71}, for $\rr,\ss$ in the simplex 
\re{simplex}
there is the following determinantal formula for
$P_\dd(\rr,\ss)$ as in \re{os2}: 
\eq
\lb{steck3}
P_\dd(\rr,\ss) = n! \, \det \left[ 
{1 (j-i+1 \ge  0) \over (j-i+1)! } \, ( s_i - r_j)^{j-i+1} _+
\right]_{1 \le i,j \le n} .
\en
The special case of \re{steck} when $s_\dd \le 1$ can be read from
\re{os}, \re{os2} and the special case of \re{steck3} with $\rr=
\vzero$ and $\ss$ the vector of partial sums of $\xx$.  The general
case of \re{steck} follows by homogeneity of $V_\dd$ from the special
case, with $x_i$ replaced by $x_i/\sigma$ for arbitrary $\sigma \ge
\sum_{i=1}^\dd x_i$.  See also Niederhausen \cite{nieder81}, where
probabilities of the form \re{steck3} are expressed in terms of
Sheffer polynomials.

\section{Sections of order cones} 
\label{ordercone}
We will obtain some results
for a class of polytopes we call ``sections of order cones''
and then show in the next section how these results apply directly to
$\Pi_\dd(\xx)$. Let $P$ be a partial ordering of the set $\{
\alpha_1,\dots, \alpha_p\}$, such that if $\alpha_i<\alpha_j$ then $i<j$. A
\emph{linear extension} of $P$ is an order-preserving bijection $\pi:P
\rightarrow [p]=\{1,2,\dots,p\}$, so if $z< z'$ in $P$ then
$\pi(z)< \pi(z')$.  We will identify $\pi$ with the
permutation (written as a word) $a_1\cdots a_p$ of $[p]$ defined by
$\pi(\alpha_{a_i}) =i$. In particular, the identity permutation
$12\cdots p$ is a linear extension of $P$. Let ${\cal L}(P)$ denote
the set of linear extensions of $P$.  Given $\pi=a_1\cdots a_p\in
{\cal L}(P)$ define ${\cal A}_\pi$ to be the set of all
order-preserving maps $f:P\rightarrow \reals$ such that
  $$ \begin{array}{c} f(\alpha_{a_1})\leq f(\alpha_{a_2}) \leq \cdots
    \leq f(\alpha_{a_p})\\[-.02em]
  f(\alpha_{a_j}) < f(\alpha_{a_{j+1}}),\ \mathrm{if}\
    a_j>a_{j+1}. \end{array} $$ 

A basic property of order-preserving maps $f:P \rightarrow\reals$ is
given by the following theorem, which is equivalent to
\cite[Lemma~4.5.3(a)]{stanley86}.

\begin{thm} {ftpp}
The set of all order-preserving maps $f:P\rightarrow \reals$ is a
disjoint union of the sets ${\cal A}_\pi$ as $\pi$ ranges over ${\cal
  L}(P)$. 
\end{thm}

For instance, if $P$ is given by Figure~\ref{fig:poset} then the
order-preserving maps $f:P\rightarrow\reals$ are partitioned by the
following seven conditions
  \beq \begin{array}{ccccccccccc}
   f(\alpha_1) & \leq & f(\alpha_2) & \leq & f(\alpha_3) & \leq &
   f(\alpha_4) & \leq & 
     f(\alpha_5) & \leq & f(\alpha_6)\nonumber\\[.1em] 
   f(\alpha_1) & \leq & f(\alpha_2) & \leq & f(\alpha_3) & \leq &
   f(\alpha_5) & < & f(\alpha_4) 
     & \leq & f(\alpha_6)\nonumber\\[.1em] 
   f(\alpha_1) & \leq & f(\alpha_3) & < & f(\alpha_2) & \leq &
   f(\alpha_4) & \leq & f(\alpha_5) 
     & \leq & f(\alpha_6)\nonumber\\[.1em] 
   f(\alpha_1) & \leq & f(\alpha_3) & < & f(\alpha_2) & \leq &
   f(\alpha_5) & < & f(\alpha_4) & 
     \leq & f(\alpha_6) \label{eq:7eq}\\[.1em] 
   f(\alpha_1) & \leq & f(\alpha_3) & \leq & f(\alpha_5) & < &
   f(\alpha_2) & \leq & f(\alpha_4) 
     & \leq & f(\alpha_6)\nonumber\\[.1em] 
   f(\alpha_2) & < & f(\alpha_1) &\leq  & f(\alpha_3) & \leq &
   f(\alpha_4) & \leq & f(\alpha_5) 
     & \leq & f(\alpha_6)\nonumber\\[.1em] 
   f(\alpha_2) & < & f(\alpha_1) & \leq & f(\alpha_3) &\leq  &
   f(\alpha_5) & < & f(\alpha_4) & 
     \leq & f(\alpha_6) \nonumber
    \end{array} \eeq

 \setlength{\unitlength}{.6pt} 
   \begin{figure}
   \begin{picture}(150,150)(-250,-10)
   \put(50,0){\circle*{5}}   
   \put(0,50){\circle*{5}}   
   \put(100,50){\circle*{5}}   
   \put(50,100){\circle*{5}}   
   \put(150,100){\circle*{5}}   
   \put(100,150){\circle*{5}}   
   \put(50,0){\line(1,1){100}}
   \put(0,50){\line(1,1){100}}
   \put(100,50){\line(-1,1){50}}
   \put(150,100){\line(-1,1){50}}
   \put(56,-5){$\alpha_1$}
   \put(106,45){$\alpha_3$}
   \put(156,95){$\alpha_5$}
   \put(-23,48){$\alpha_2$}
   \put(24,98){$\alpha_4$}
   \put(74,148){$\alpha_6$}
\end{picture}
\caption{A partially ordered set}
\label{fig:poset}
\end{figure}
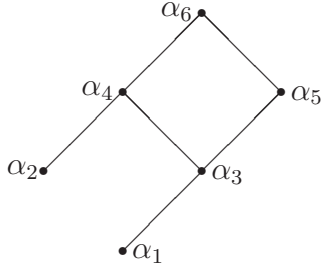

Define the \emph{order cone} ${\cal C}(P)$ of the poset $P$ to be the
set of all order-preserving maps $f:P\rightarrow \nreals$. Thus ${\cal
  C}(P)$ is a pointed polyhedral cone in the space $\reals^P$.  Assume
now that $P$ has a unique maximal element $\hat{1}$, and let
$t_1<\cdots<t_\dd=\hat{1}$ be a chain $C$ in $P$. (With a little more
work we could relax the assumption that $C$ is a chain. The condition
that $t_\dd=\hat{1}$ entails no real loss of generality since we can
just adjoin a $\hat{1}$ to $P$ and include it in $C$.)  Let
$x_1,\dots,x_\dd$ be nonnegative 
real numbers. Set $u_i=x_1+\cdots+x_i$ and $\uu =(u_1,\dots, u_\dd)$.  Let
$W_{\uu}$ denote the subspace of $\reals^P$ defined by $f(t_i)=u_i$ for
$1\leq i\leq\dd$.  Define the \emph{order cone section} ${\cal
  C}_C(P,\uu)$ to be the intersection ${\cal C}(P)\cap W_\suu$, restricted
to the coordinates $P-C$. (The restriction to the coordinates $P-C$
merely deletes constant coordinates and has no effect on the geometric
and combinatorial structure of ${\cal C}(P)\cap W_\suu$.) Equivalently,
${\cal C}_C(P,\uu)$ is the set of all order-preserving maps $f:P-C
\rightarrow\nreals$ such that the extension of $f$ to $P$ defined by
$f(t_i)=u_i$ remains order-preserving. Note that ${\cal C}_C(P,\uu)$
is bounded since for all $s\in P-C$ and all $f\in{\cal C}_C(P,\uu)$ we
have $0\leq f(s)\leq u_\dd$. Thus ${\cal C}_C(P,\uu)$ is a
convex polytope contained in $\reals^{P-C}$. Moreover, $\dim {\cal
  C}_C(P,\uu)=|P-C|$ provided each $x_i>0$ (or in certain other
situations, such as when no element of $P-C$ is greater than $t_1$).

There is an alternative way to view the polytope ${\cal C}_C(P,\uu)$. Let
${\cal P}_1,\dots,{\cal P}_\dd$ be convex polytopes (or just convex
bodies) in the same ambient space $\reals^m$, and let $x_1,\dots,x_\dd
\in\nreals$. Define the \emph{Minkowski sum} (or more accurately,
\emph{Minkowski linear combination})
  $$ x_1{\cal P}_1+\cdots+x_\dd{\cal P}_\dd = \{ x_1 X_1+\cdots+
      x_\dd X_\dd\,:\,X_i\in {\cal P}_i\}. $$
Then ${\cal Q}=x_1{\cal P}_1+\cdots+x_\dd{\cal P}_\dd$ is a convex
polytope that was first investigated by Minkowski (at least for $m\leq
3$) and whose study belongs to the subject of \emph{integral geometry}
(e.g., \cite{schneider}). In particular, the $m$-dimensional volume of
${\cal Q}$ has the form
  $$ \mathrm{Vol}({\cal Q}) = \sum_{{a_1+\cdots+a_\dd=m\atop a_i\in
      \mathbb{N}}} {m\choose a_1,\dots,a_\dd}V({\cal P}_1^{a_1},
      \dots,{\cal P}_\dd^{a_\dd})x_1^{a_1}\cdots x_\dd^{a_\dd}, $$
where $V({\cal P}_1^{a_1},\dots,{\cal P}_\dd^{a_\dd})\in\nreals$. 
These numbers are known as the \emph{mixed volumes} of the polytopes
${\cal P}_1,\dots,{\cal P}_\dd$ and have been extensively
investigated. 

Now suppose that ${\cal P}_1,\dots,{\cal P}_\dd$ are
\emph{integer polytopes} (i.e., their vertices have integer
coordinates) in $\reals^m$, and let
$x_1,\dots,x_\dd\in\mathbb{N}$. Given any integer polytope ${\cal
  P}\subset \reals^m$,
write
      $$ N({\cal P}) = \#({\cal P}\cap \mathbb{Z}^m), $$
the number of integer points in ${\cal P}$. Then we call $N(x_1{\cal
  P}_1 +\cdots+x_\dd{\cal P}_\dd)$, regarded as a function of
$x_1,\dots,x_\dd\in\mathbb{N}$, the \emph{mixed lattice point
  enumerator} of ${\cal P}_1,\dots,{\cal P}_\dd$. It was shown by
McMullen \cite{mcmullen} (see also \cite{mcmullen2}\cite{mcmullen3}
for two related survey articles) that $N(x_1{\cal P}_1
+\cdots+x_\dd{\cal P}_\dd)$ is a polynomial in $x_1,\dots,x_\dd$ (with
rational coefficients) of total degree at most $m$. Moreover, the
terms of degree $m$ are given by Vol$(x_1{\cal P}_1
+\cdots+x_\dd{\cal P}_\dd)$. Hence the coefficients of the terms of
degree $m$ are nonnegative, but in general the coefficients of
$N(x_1{\cal P}_1 +\cdots+x_\dd{\cal P}_\dd)$ may be negative. In the
special case $\dd=1$, the mixed lattice point enumerator $N(x{\cal
  P})$ is called the \emph{Ehrhart polynomial} of the integer polytope
${\cal P}$ and is denoted $i({\cal P},x)$. An introduction to Ehrhart
polynomials appears in \cite[pp.\ 235--241]{stanley86}. 

Define the \emph{order polytope} ${\cal O}(P)$ of the finite poset $P$
to be the set of all order-preserving maps $f:P\rightarrow [0,1] = \{
x\in\reals\,:\, 0\leq x\leq 1\}$. Thus ${\cal O}(P)$ is a convex
polytope in $\reals^P$ of dimension $|P|$. The basic properties of
order polytopes are developed in \cite{stanley2pp}.

\begin{thm} {popms}
Given $P$, $C$, and $\uu$ as above, so $u_i=x_1+\cdots+x_i$, let 
  $$ P_i = \{ s\in P-C\,:\,s\not <t_{i-1} \} $$
(with $P_1=P-C$). Regard the order polytope ${\cal O}(P_i)$ as lying in
$\reals^{P-C}$ by setting coordinates indexed by elements of
$(P-C)-P_i$ equal to $0$. Then
  $$ {\cal C}_C(P,\uu) = x_1{\cal O}(P_1)+x_2{\cal O}(P_2)+\cdots +
       x_\dd{\cal O}(P_\dd). $$
\end{thm}

\proof
We can regard ${\cal O}(P_i)$ as the set of order preserving maps
$f:P-C\rightarrow [0,1]$ such that $f(s)=0$ if $s<t_{i-1}$. From this
it is clear that every element of $x_1{\cal O}(P_1)+x_2{\cal
O}(P_2)+\cdots + x_\dd{\cal O}(P_\dd)$ is an order-preserving map
$g:P-C\rightarrow \nreals$ such that the extension of $g$ to $P$
defined by $g(t_i)=x_1+\cdots+x_i$ remains order-preserving. Hence 
  $$ {\cal C}_C(P,\uu) \supseteq x_1{\cal O}(P_1)+x_2{\cal
       O}(P_2)+\cdots +  x_\dd{\cal O}(P_\dd). $$
For the converse, we may assume (by deleting elements of $P$ if
necessary) that each $x_i>0$. let $f\in{\cal C}_C(P,\uu)$. Let $s\in
P_C$ and define $g_1(s)=f(s)$ and $f_1(s)=\min(1,x_1^{-1}g_1(s))$. Set
  $$ g_2(s)=g_1(s)-x_1f_1(s)=\max(g_1(s)-x_1,0). $$
Now let $f_2(s)=\min(1,x_2^{-1}g_2(s))$ and set
  $$ g_3(s)=g_2(s)-x_2f_2(s)=\max(g_2(s)-x_2,0). $$
Continuing in this way gives functions $f_1,f_2,\dots,f_\dd$, for which
it can be checked that $f_i\in{\cal O}(P_i)$ and 
  $$ f=x_1f_1+\cdots+x_\dd f_\dd, $$
so 
  $$ {\cal C}_C(P,\uu) \subseteq x_1{\cal O}(P_1)+x_2{\cal
       O}(P_2)+\cdots +  x_\dd{\cal O}(P_\dd). $$
\endpf

We now want to give a formula for the number of integer points in
${\cal C}_C(P,\uu)$, which by Theorem~\ref{popms} is just the mixed
lattice point enumerator of the polytopes ${\cal O}(P_i)$. Let $C$ be
the chain $t_1<\cdots<t_\dd=\hat{1}$ as above. Given
$\pi=a_1 \cdots a_p\in {\cal L}(P)$, write $h_i(\pi)$ for the
\emph{height} of $t_i$ in $\pi$, i.e.,
$t_i=\pi^{-1}(a_{h_i(\pi)})$. Thus $1\leq 
h_1(\pi)<\cdots<h_\dd(\pi)=p$. Also write 
  $$ d_i(\pi) = \#\{ j\,:\, h_{i-1}(\pi)\leq j<h_i(\pi),\ a_j>a_{j+1}
      \}, $$
where we set $h_0(\pi)=0$ and $a_0=0$. Thus $d_i(\pi)$ is the number
of \emph{descents} of $\pi$ appearing between $h_{i-1}(\pi)$ and
$h_i(\pi)$. Recall (e.g., \cite[{\S}1.2]{stanley86}) that
the number of ways to choose $j$ objects with repetition from a set of
$k$ objects is given by
  \beq \bbc{k}{j}={k+j-1\choose j}=\frac{k(k+1)\cdots(k+j-1)}{j!}.
    \label{eq:bbc} \eeq
Regarding
$\left(\hspace{-.4em}\left({k\atop j}\right)\hspace{-.4em}\right)$ as a
  polynomial in $k\in\mathbb{Z}$, note that 
$\left(\hspace{-.4em}\left({k\atop j}\right)\hspace{-.4em}\right)=0$ for
  $-j+1\leq k\leq 0$. 

\begin{thm} {propmle}
We have
  \beq N({\cal C}_C(P,\uu)) = \sum_{\pi\in{\cal L}(P)} \prod_{i=1}^{\dd-1}
   \bbc{x_i-d_i(\pi)+1}{h_i(\pi)-h_{i-1}(\pi)-1}. \label{eq:noc}
  \eeq
\end{thm}

\proof
Fix $\pi=a_1\cdots a_p\in{\cal L}(P)$. Write $h_i=h_i(\pi)$ and
$d_i=d_i(\pi)$. Let $f:P\rightarrow\reals$ be an order-preserving map
such that (a) $f\in {\cal A}_\pi$, (b) $f(t_i)=u_i=x_1+\cdots+x_i$,
and (c) the restriction $f|_{P-C}$ of $f$ to $P-C$ satisfies $f|_{P-C}
\in {\cal C}_C(P,\uu)$. If we write $c_i=f(\alpha_{a_i})$, then for fixed
$\pi$ it follows from Theorem~\ref{ftpp} that the integer points
$f|_{P-C}\in {\cal C}_C(P,\uu)$, where $f$ satisfies (a) and (b), are
given by
  $$ 0\leq c_1\leq c_2\leq\cdots\leq c_{h_1}=x_1\leq c_{h_1+1}\leq
     \cdots \leq c_{h_2}=x_1+x_2 $$
\vspace{-.3in}
  \beq \qquad \leq \cdots\leq c_p=x_1+\cdots+x_\dd
  \label{eq:abc} \eeq
\vspace{-.05in}
  \beq  c_j<c_{j+1}\ \mbox{if}\ a_j>a_{j+1}. \label{eq:st} \eeq
Let $\alpha, \beta,m\in\mathbb{N}$ and $0\leq j_1<j_2<\cdots<j_q\leq
m$. Elementary combinatorial reasoning shows that the number of
integer vectors $(r_1,\dots,r_m)$ satisfying
  $$ \alpha=r_0\leq r_1\leq \cdots\leq r_m\leq r_{m+1}=\alpha+
      \beta $$ 
\vspace{-.3em}
    $$ r_{j_i}<r_{j_i}+1\ \mbox{for}\ 1\leq i\leq q $$
is equal to $\left(\hspace{-.3em}\left({\beta-q+1\atop
      m}\right)\hspace{-.3em}\right)$. Hence the number of integer 
sequences satisfying (\ref{eq:abc}) and (\ref{eq:st}) is given by
  $$ \bbc{x_1-d_1+1}{h_1-1}\bbc{x_2-d_2+1}{h_2-h_1-1}\cdots
      \bbc{x_\dd-d_\dd+1}{h_\dd-h_{\dd-1}-1}. $$
Summing over all $\pi\in {\cal L}(P)$ yields (\ref{eq:noc}).
\endpf

\begin{exm}{7linex}
Let $P$ be given by Figure~\ref{fig:poset}, and let $t_1=\alpha_1$,
$t_2=\alpha_3$, and $t_3=\alpha_6$. The conditions in equation (\ref{eq:7eq})
become in the notation of the above proof as follows:
  $$ 0\leq c_1 =x_1\leq c_2\leq  c_3=x_1+x_2\leq  c_4 \leq c_5 \leq
  c_6 = x_1 +x_2+x_3 $$ 
\vspace{-.1in}
  $$ 0\leq c_1=x_1\leq  c_2\leq  c_3=x_1+x_2\leq  c_4< c_5 \leq c_6 =
  x_1 +x_2+x_3 $$ 
\vspace{-.1in}
  $$ 0\leq c_1=x_1\leq  c_2=x_1+x_2< c_3\leq  c_4\leq  c_5 \leq c_6 =
  x_1 +x_2+x_3 $$
\vspace{-.1in} 
  $$ 0\leq c_1=x_1\leq  c_2=x_1+x_2< c_3\leq  c_4 <c_5 \leq c_6 = x_1
  +x_2+x_3 $$ 
\vspace{-.1in}
  $$ 0\leq c_1=x_1\leq  c_2=x_1+x_2\leq  c_3< c_4\leq  c_5 \leq c_6 =
  x_1 +x_2+x_3 $$ 
\vspace{-.1in}
  $$ 0\leq c_1< c_2=x_1\leq  c_3=x_1+x_2\leq  c_4\leq  c_5 \leq c_6 =
  x_1 +x_2+x_3 $$ 
\vspace{-.1in}
  $$ 0\leq c_1 <c_2=x_1\leq  c_3=x_1+x_2\leq  c_4< c_5 \leq c_6 = x_1
  +x_2+x_3, $$ 
yielding
  $$ N({\cal C}_C(P,\uu))=\bbc{x_2+1}{1}\bbc{x_3+1}{2}+
   \bbc{x_2+1}{1}\bbc{x_3}{2}+\bbc{x_3}{3} $$
\vspace{-.1in}
   $$ \quad +\bbc{x_3-1}{3}+
   \bbc{x_3}{3}+\bbc{x_1}{1}\bbc{x_3+1}{2}+\bbc{x_1}{1}\bbc{x_3}{2}. $$
\end{exm}

We mentioned earlier that the terms of highest degree (here of degree
$|P-C|$) of $N(x_1{\cal P}_1 +\cdots+x_\dd{\cal P}_\dd)$ are given by
Vol$(x_1{\cal P}_1+\cdots+x_\dd{\cal P}_\dd)$. Hence we obtain from
Theorem~\ref{propmle} the following result.

\begin{crl}{projvol}
The volume of ${\cal C}_C(P,\uu)$ is given by
  \beq \mathrm{Vol}({\cal C}_C(P,\uu)) =  \sum_{\pi\in{\cal L}(P)}
         \prod_{i=1}^\dd \frac{x_i^{h_i(\pi)-h_{i-1}(\pi)}}
        {(h_i(\pi)-h_{i-1}(\pi))!}. \label{eq:volo} \eeq
Thus if $m=|P-C|$ then the mixed volume $m!\cdot V({\cal
O}(P_1)^{a_1},\dots, {\cal O}(P_\dd)^{a_\dd})$ is equal to the number of
linear extensions $\pi\in{\cal L}(P)$ such that $t_i$ has height
$a_1+\cdots +a_i$ in $\pi$, for $1\leq i\leq \dd$.
\end{crl}

The case $\dd=2$ of Corollary~\ref{projvol} (or equivalently the case
$\dd=1$ where $t_1$ can be any element of $P$, not just the top element)
appears in \cite[(16)]{stanley2pp}. 

The \emph{product} of two polytopes ${\cal P}\in \reals^p$ and ${\cal
Q}\in\reals^q$ is defined to be their cartesian product ${\cal
P}\times {\cal Q}\in\reals^{p+q}$. If $\bar{L}({\cal P})$ denotes
the poset of \emph{nonempty} faces of ${\cal P}$, then $\bar{L}({\cal
P}\times {\cal Q})=\bar{L}({\cal P})\times \bar{L}({\cal Q})$ (see
Ziegler \cite[pp.\ 9--10]{ziegler95}). If ${\cal P}$ is a $d$-simplex,
then $\bar{L}({\cal P})$ is just a boolean algebra of rank $d$ with
the minimum element removed. Moreover, the product of $\dd$
one-dimensional simplices is combinatorially equivalent (even affinely
equivalent) to a $d$-cube.
If $\pi=a_1\cdots a_p\in{\cal L}_P$, then define $\Lambda_\pi$ to be
the subset of ${\cal C}_C(P,\uu)$ given by equation (\ref{eq:abc}). Thus
when each $x_i>0$ we have that $\Lambda_\pi$ is a product of simplices
of dimensions $h_1-1$, $h_2-h_1-1,\dots,h_p-h_1-1$, and 
  $$ \mathrm{Vol}(\Lambda_\pi)= \prod_{i=1}^\dd
    \frac{x_i^{h_i(\pi)-h_{i-1}(\pi)}}{(h_i(\pi)-h_{i-1}(\pi))!}. $$
Moreover, the $\Lambda_\pi$'s form the chambers of a polyhedral
decomposition $\Omega_C(P,\uu)$ of ${\cal C}_C(P,\uu)$. We regard
$\Omega_C(P,\uu)$ as the set of all faces of the $\Lambda_\pi$'s
(including the $\Lambda_\pi$'s themselves), partially ordered by
inclusion. Note that the formula (\ref{eq:volo}) corresponds to an
explicit decomposition of ${\cal C}_C(P,\uu)$ into ``nice'' pieces
(products of simplices) whose volumes are the terms in (\ref{eq:volo}).

Our next result concerns the 
combinatorial structure of the decomposition of ${\cal C}_C(P,\uu)$ into
the chambers $\Lambda_\pi$. First we review some information from
\cite[{\S}5]{stanley2pp} about the cone ${\cal C}(P)$ of all
order-preserving maps $f:P\rightarrow \nreals$. (The paper
\cite{stanley2pp} actually deals with the order complex ${\cal O}(P)$
rather than the cone ${\cal C}(P)$, but this does not affect our
arguments.) Recall (e.g., \cite[p.\ 100]{stanley86}) that an
\emph{order ideal} $I$ of $P$ is a subset of $P$ such that if $t\in I$
and $s<t$, then $s\in I$. The poset (actually a distributive lattice)
of all order ideals of $P$, ordered by inclusion, is denoted
$J(P)$. Given a chain $K:\emptyset=I_0< I_1<\cdots<I_k=P$ in $J(P)$,
define ${\cal C}_K(P)$ to consist of all $f:P\rightarrow\nreals$
satisfying 
 \beq 0\leq f(I_1)\leq f(I_2-I_1)\leq \cdots \leq f(I_k-I_{k-1}), 
   \label{eq:jpchain} \eeq
where $f(S)$ denotes the common value of $f$ at all
the elements of the subset $S$ of $P$. Clearly ${\cal C}_K(P)$ is a
$k$-dimensional cone in $\reals^P$. It is not hard to see that the set
$\Omega(P) =\{ {\cal C}_K(P)\,:\,K\mbox{ is a chain in }J(P)\mbox{
  containing }\emptyset\mbox{ and }P\}$ is a triangulation of ${\cal
  C}(P)$. The chambers (maximal faces) of $\Omega(P)$ consist of the
cones
   $$ 0\leq f(\alpha_{a_1})\leq \cdots\leq f(\alpha_{a_p}), $$
where $\pi=a_1\cdots a_p\in {\cal L}(P)$. Moreover, ${\cal C}_K(P)$ is
an \emph{interior} face of $\Omega(P)$ (i.e., does not lie on the
boundary) if and only if each subset $I_i-I_{i-1}$ of equation
(\ref{eq:jpchain}) is an \emph{antichain}, i.e., no two distinct
elements of $I_i-I_{i-1}$ are comparable. Such chains of $J(P)$ are
called \emph{Loewy chains}. Let $\Omega^{\circ}(P)$
denote the set of interior faces of $\Omega(P)$ regarded as a
partially ordered set under inclusion. Thus $\Omega^\circ(P)$ is
isomorphic to the set of Loewy chains of $J(P)$, ordered by
inclusion. Similarly, we let $\Omega_C^\circ(P,\uu)$ denote the set of
interior faces of the polyhedral decomposition $\Omega_C(P,\uu)$. 

\begin{thm}{thm:decomp}
  Let $W_\suu$ denote the subspace of $\reals^P$ given by $f(t_i)=u_i$,
  $1\leq i\leq\dd$. Define a map $\phi:\Omega^\circ(P) \rightarrow
  \Omega_C^\circ(P,\uu)$ by letting $\phi({\cal C}_K(P))$ equal
  $\phi{\cal}_K(P)\cap W_\suu$ restricted to the coordinates $P-C$. Then
  $\phi$ is an isomorphism of posets.
\end{thm}

\proof 
Let (\ref{eq:jpchain}) define an interior face ${\cal C}_K(P)$ of
${\cal C}(P)$, so $\emptyset=I_0< I_1<\cdots<I_k=P$ is a Loewy
chain. Thus each set $I_j-I_{j-1}$ contains at most one element of the
chain $C:\,t_1<\cdots<t_\dd$. Let $t_i\in I_{j_i}-I_{j_i-1}$. (In
particular, $j_\dd=k$ since $t_\dd=\hat{1}$.) Then
$\phi({\cal C}_K(P))$ is defined by the equations
    $$ \begin{array}{c} 0\leq f(I_1)\leq f(I_2-I_1)\leq \cdots \leq
    f(I_{j_1}-I_{j_1-1})=u_1\\[-.02em]
   \ \ \leq f(I_{j_1+1}-I_{j_1})\leq \cdots \leq
    f(I_{j_2}-I_{j_2-1})=u_2\leq \cdots\leq f(I_k-I_{k-1})=u_\dd. 
    \end{array} $$
It follows immediately that $\phi$ is a bijection, and that two Loewy
chains $K$ and $K'$ satisfy $K\subseteq K'$ if and only if 
$\phi({\cal C}_K(P))\subseteq \phi({\cal C}_{K'}(P))$. Hence $\phi$ is
a poset isomorphism. 
\endpf

The point of Theorem~\ref{thm:decomp} is that it gives a simple
combinatorial description (namely, the poset $\Omega^\circ(P)$, which
is isomorphic to the set of Loewy chains of $J(P)$ under inclusion) of
the geometrically defined poset $\Omega_C^\circ(P,\uu)$. Note that 
$\Omega^\circ(P)$ depends only on $P$, not on the chain $C$.

\section{$\Pi_\dd(\xx)$ as a section of an order cone}
\label{pinsec}
In this section we will apply the theory developed in the previous
section to $\Pi_\dd(\xx)$. Let us say that two integer polytopes
${\cal P}\subset \reals^k$ and ${\cal Q}\subset \reals^m$ are
\emph{integrally equivalent} if there is an affine transformation
$\varphi:\reals^k\rightarrow \reals^m$ whose restriction to ${\cal P}$
is a bijection $\varphi:{\cal P}\rightarrow {\cal Q}$, and such that
if aff denotes affine span, then $\varphi$ restricted to
$\mathbb{Z}^k\cap \mathrm{aff}({\cal P})$ is a bijection
$\varphi:\mathbb{Z}^k\cap \mathrm{aff}({\cal P}) \rightarrow
\mathbb{Z}^m\cap \mathrm{aff}({\cal Q})$. It follows that ${\cal P}$
and ${\cal Q}$ have the same combinatorial type and the same
``integral structure,'' and hence the same volume, Ehrhart polynomial,
etc.

 Now let $\mbox{\boldmath$i$}$ denote an $i$-element chain, and
let $Q_\dd={\mbox{\boldmath$2$}}\times{\mbox{\boldmath$\dd$}}$, the
product of a two-element 
chain with an $\dd$-element chain. We regard the elements of $Q_\dd$ as
$\alpha_1,\dots,\alpha_{2\dd}$ with $\alpha_1<\cdots<\alpha_\dd$,
$\alpha_{\dd+1}<\cdots <
\alpha_{2\dd}$, and $\alpha_i<\alpha_{\dd+i}$ for $1\leq i\leq \dd$. Let
$t_i=\alpha_{\dd+i}$, and let $C$ be the chain $t_1<\cdots<t_\dd$. As in
the previous section let 
$x_1,\dots,x_\dd\geq 0$, and set $u_i=x_1+\cdots+x_i$. The polytope
${\cal C}_C(Q_\dd,\uu)\subset \reals^{Q_\dd-C}\cong\reals^\dd$ thus by
definition is given by the equations 
   $$ 0\leq f_1\leq \cdots\leq f_\dd $$
   $$ f_i\leq u_i,\ 1\leq i\leq \dd. $$
Let $y_i=f_i-f_{i-1}$ (with $f_0=0$). Then the above equations become
   $$ y_i\geq 0,\ 1\leq i\leq \dd $$
   $$ y_1+\cdots+y_i\leq x_1+\cdots+x_\dd. $$ 
These are just the equations for $\Pi_\dd(\xx)$. The
transformation $y_i=f_i-f_{i-1}$ induces an integral equivalence
between ${\cal C}_C(Q_\dd,\uu)$ and $\Pi_\dd(\uu)$. Hence the results of
the above section, when specialized to $P=Q_\dd$, are directly
applicable to $\Pi_\dd(\xx)$.

Theorem~\ref{popms} expresses ${\cal C}_C(P,\uu)$ as a Minkowski
linear combination of order polytopes ${\cal O}(P_i)$. In the present
situation, where $P=\bm{2}\times\bm{\dd}$, the poset $P_i$ is just the
chain $\alpha_i<\alpha_{i+1}<\cdots<\alpha_\dd$. The order polytope
${\cal O}(P_i)$ is defined by the conditions
  $$ f_1=\cdots=f_{i-1}=0,\ \  0\leq f_i\leq \cdots\leq f_\dd\leq 1. $$
This is just a simplex of dimension $\dd-i+1$ with vertices
$(0^j,1^{\dd-j})$, $i-1\leq j\leq \dd$, where $(0^j,1^{\dd-j})$ denotes a
vector of $j$ 0's followed by $\dd-j$ 1's. Switching to the $y$
coordinates (i.e., $y_i=f_i-f_{i-1}$) yields the following result.

\begin{thm} {pinms}
Let $\tau_i$ be the $(\dd-i+1)$-dimensional simplex in $\reals^\dd$
defined by
  $$ \begin{array}{c} y_1=\cdots=y_{i-1}=0\\ 
      y_i\geq 0,\dots, y_\dd\geq 0\\
      y_i+\cdots+y_\dd\leq 1, \end{array} $$
with vertices $(0^{j-1},1,0^{\dd-j})$ for $i\leq j\leq \dd$, and
$(0,0,\dots,0)$. Then
  $$ \Pi_\dd(\xx) = x_1\tau_1+x_2\tau_2+\cdots+x_\dd\tau_\dd. $$
\end{thm}

Consider the set ${\cal L}(Q_\dd)$ of linear extensions of $Q_\dd$. A
linear extension $\pi=a_1\dots a_{2\dd}\in {\cal L}(Q_\dd)$ is
uniquely determined by the positions of $\dd+1,\dots,2\dd$ (since
$1,\dots,\dd$ must appear in increasing order). If $a_{j_i}=\dd+i$ for
$1\leq i\leq\dd$, then $1\leq j_1<\cdots<j_\dd=2\dd$ and $j_i\geq
2i$. The number of such sequences is just the Catalan number $C_\dd =
\frac{1}{\dd+1}{2\dd \choose \dd}$ (see e.g.\ 
\cite[Exercise~6.19(t)]{stanleyv2}, which is a minor variation). If we
set $k_i=j_i-j_{i-1}$ (with $j_0=0$), then the sequences
$\bm{k}=(k_1,\dots,k_\dd)$ are just those of equation (\ref{kdef}).
Moreover, in the linear extension $a_1\cdots a_{2\dd}$ there are no
descents to the left of $\dd+1$, and there is exactly one descent
between $\dd+i$ and $\dd+i+1$ provided that $k_{i+1}-k_i\geq 2$. (If
$k_{i+1}-k_i=1$ then there are no descents between $\dd+i$ and
$\dd+i+1$.) By Theorem~\ref{propmle} we conclude
  \beq N(\Pi_\dd(\xx)) = \sum_{\skk\in K_\dd} \bbc{x_1+1}{k_1}
    \prod_{i=2}^\dd \bbc{x_i}{k_i}, \label{Npin} \eeq
where $K_\dd$ is given by (\ref{kdef}).
Taking terms of highest degree yields Theorem~\ref{volume}. Thus we 
have obtained an explicit decomposition of $\Pi_\dd(\xx)$ into products
of simplices whose volumes are the terms in (\ref{vform}). (A
completely different such decomposition will be given in
Section~\ref{sec:assoc}.) 
Moreover, Theorem~\ref{thm:decomp} gives the combinatorial structure
of the interior faces of this decomposition.

\textsc{Note.} Equation~(\ref{Npin}) was obtained independently by Ira
Gessel (private communication) by a different method.

 \setlength{\unitlength}{.6pt} 
   \begin{figure}
   \begin{picture}(150,150)(-250,-10)
   \put(50,0){\circle*{5}}   
   \put(0,50){\circle*{5}}   
   \put(100,50){\circle*{5}}   
   \put(50,100){\circle*{5}}   
   \put(150,100){\circle*{5}}   
   \put(100,150){\circle*{5}}   
   \put(50,0){\line(-1,1){50}}
   \put(50,0){\line(1,1){100}}
   \put(0,50){\line(1,1){100}}
   \put(100,50){\line(-1,1){50}}
   \put(150,100){\line(-1,1){50}}
   \put(56,-5){$\alpha_1$}
   \put(106,45){$\alpha_2$}
   \put(156,95){$\alpha_3$}
   \put(-25,48){$\alpha_4$}
   \put(24,98){$\alpha_5$}
   \put(74,148){$\alpha_6$}
\end{picture}
\caption{The poset $Q_3={\mbox{\boldmath$2$}}\times{\mbox{\boldmath$3$}}$}
\label{fig:p3}
\end{figure}
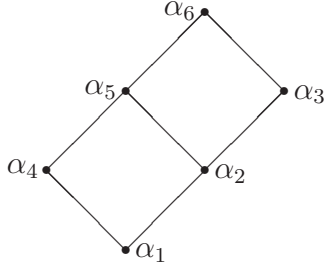

Let us illustrate the above discussion with the case $\dd=3$. The poset
$Q_3$ is shown in Figure~\ref{fig:p3}. The linear extensions of $Q_3$
are given as follows, with the elements $4,5,6$ corresponding to the
chain $C$ shown in boldface:
  $$ \begin{array}{l}
    123{\mbox{\boldmath$4$}}{\mbox{\boldmath$5$}}{\mbox{\boldmath$6$}}\\
    12{\mbox{\boldmath$4$}}3{\mbox{\boldmath$5$}}{\mbox{\boldmath$6$}}\\
    12{\mbox{\boldmath$4$}}{\mbox{\boldmath$5$}}3{\mbox{\boldmath$6$}}\\
    1{\mbox{\boldmath$4$}}23{\mbox{\boldmath$5$}}{\mbox{\boldmath$6$}}\\
    1{\mbox{\boldmath$4$}}2{\mbox{\boldmath$5$}}3{\mbox{\boldmath$6$}}
  \end{array} $$
Hence the points $(y_1,y_2,y_3)\in\Pi_3(\xx)$ are decomposed into the
sets 
  \beq \begin{array}{c} 0\leq y_1\leq y_2\leq y_3\leq x_1\\
   0\leq y_1\leq y_2\leq x_1<y_3\leq x_1+x_2\\
   0\leq y_1\leq y_2\leq x_1\leq x_1+x_2<y_3\leq x_1+x_2+x_3\\
   0\leq y_1\leq x_1<y_2\leq y_3\leq x_1+x_2\\
   0\leq y_1\leq x_1<y_2\leq x_1+x_2<y_3\leq x_1+x_2+x_3,
  \end{array} \label{eq:exdecp3} \eeq
yielding
  $$ N(\Pi_3(\xx)) = \bbc{x_1+1}{3}+\bbc{x_1+1}{2}\bbc{x_2}{1}
                   +\bbc{x_1+1}{2}\bbc{x_3}{1} $$
\vspace{-.1in}
  $$ \qquad + \bbc{x_1+1}{1}\bbc{x_2}{2}+\bbc{x_1+1}{1}
      \bbc{x_2}{1}\bbc{x_3}{1}. $$
Theorem~\ref{thm:decomp} allows us to describe the incidence relations
among the faces of the decomposition of $\Pi_3(\xx)$ whose chambers are
the closures of the five sets in equation (\ref{eq:exdecp3}). The
lattice $J(Q_3)$ of order ideals of $Q_3$ has five maximal
chains. This lattice is shown in Figure~\ref{fig:oip3}, with elements
labeled $a,b,\dots, j$. The elements $a,b,i,j$ appear in every Loewy
chain of $J(Q_3)$ and can be ignored. The simplicial complex of chains
of $J(P)$ (with $a,b,i,j$ removed) is shown in
Figure~\ref{fig:p3decomp}(a). The Loewy chains correspond to the
interior faces, of which five have dimension 2, five have dimension 1,
and one has dimension 0. Figure~\ref{fig:p3decomp} shows the ``dual
complex'' of the interior faces. This gives the incidence relations
among the five chambers of the decomposition of $\Pi_3(\xx)$ into five
products of simplices obtained from $\Omega_C^\circ(P,\uu)$ by the
change of coordinates $y_i=f_i-f_{i-1}$ discussed above. For a
picture, see the second subdivision of $\Pi_3(\xx)$ in
Figure~\ref{fig3}. 

 \setlength{\unitlength}{.6pt} 
   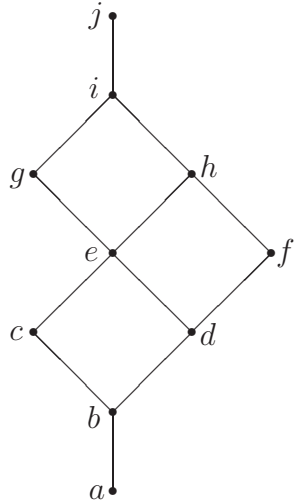
\begin{figure}
   \begin{picture}(150,320)(-250,-10)
  \put(50,0){\circle*{5}}
  \put(50,50){\circle*{5}}
  \put(0,100){\circle*{5}}
  \put(100,100){\circle*{5}}
  \put(50,150){\circle*{5}}
  \put(150,150){\circle*{5}}
  \put(0,200){\circle*{5}}
  \put(100,200){\circle*{5}}
  \put(50,250){\circle*{5}}
  \put(50,300){\circle*{5}}
  \put(50,0){\line(0,1){50}}
  \put(50,50){\line(1,1){100}}
  \put(50,50){\line(-1,1){50}}
  \put(0,100){\line(1,1){100}}
  \put(100,100){\line(-1,1){100}}
  \put(150,150){\line(-1,1){100}}
  \put(0,200){\line(1,1){50}}
  \put(50,250){\line(0,1){50}}
  \put(35,-5){$a$}
  \put(34,40){$b$}
  \put(-15,95){$c$}
  \put(105,92){$d$}
  \put(32,145){$e$}
  \put(153,145){$f$}
  \put(-15,195){$g$}
  \put(105,198){$h$}
  \put(35,248){$i$}
  \put(35,295){$j$}
\end{picture}
\caption{The lattice $J(Q_3)$ of order ideals of $Q_3$}
\label{fig:oip3}
\end{figure}

\begin{figure}
\centerline{\epsfig{figure=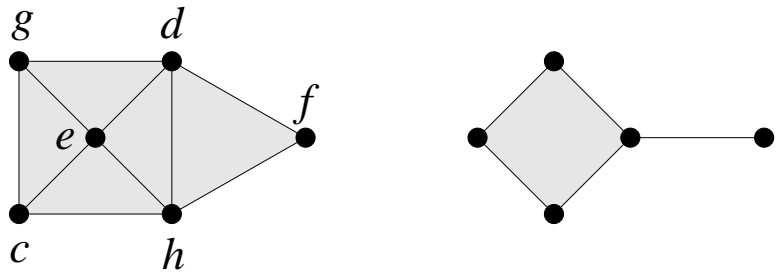}}
\caption{The order complex of $J(Q_3)$ with $a,b,i,j$ omitted,  and the
interior face dual complex}
\label{fig:p3decomp}
\end{figure}

We mentioned earlier that in general the coefficients of the mixed
lattice point enumerator $N(x_1{\cal P}_1 +\cdots+x_\dd{\cal P}_\dd)$
may be negative. The polytope $\Pi_\dd(\xx)$ is an exception, however,
and in fact satisfies a slightly stronger property.

\begin{crl}{nonneg}
The polynomial $N(\Pi_\dd(x_1-1,x_2,\dots,x_\dd))$ has nonnegative
coefficients.
\end{crl}

\proof
Immediate from equation (\ref{Npin}), since the polynomial
$\left(\hspace{-.3em}\left({t\atop i}\right)\hspace{-.3em}\right)$ has
nonnegative coefficients.
\endpf

\textsc{Note.} One can also think of ${\cal C}_C(Q_\dd,\uu)$ as the
``polytope of fractional shapes contained in the shape
$(u_\dd,u_{\dd-1},\dots,u_1)$.'' In general, let
$\lambda=(\lambda_1,\dots,\lambda_\dd)$ be a partition, i.e.,
$\lambda_i\in\mathbb{N}$ and $\lambda_1\geq \cdots\geq \lambda_\dd$,
which we also call a \emph{shape}. We say that a shape $\mu=(\mu_1,
\dots,\mu_\dd)$ is \emph{contained in} $\lambda$ if $\mu_i\leq
\lambda_i$ for all $i$. (This partial ordering on shapes defines
\emph{Young's lattice} \cite[Exer.\ 3.63]{stanley86}. Additional
properties of Young's lattice may be found in various places in
\cite{stanleyv2}.) If we relax the conditions that the $\lambda_i$'s
are integers but only require them to be real
(with $\lambda_1\geq \cdots\geq \lambda_\dd\geq 0$), then we can 
think of $\lambda$ as a ``fractional shape.'' Thus ${\cal
C}_C(Q_\dd,\uu)$ just consists of the fractional shapes contained in
the shape $(u_\dd,u_{\dd-1},\dots,u_1)$.

\section{Connections with parking functions and plane partitions.}
\label{sec:pfpp}
There are two additional interpretations of the volume and lattice
point enumerator of $\Pi_\dd(\xx)$ that we wish to discuss. The first
concerns the subject of parking functions, originally defined by
Konheim and Weiss \cite{kon-weiss}.  A \emph{parking function} of
length $\dd$ may be defined as a sequence $(a_1,\dots,a_\dd)$ of
positive integers whose increasing rearrangement $b_1\leq\cdots\leq
b_\dd$ satisfies $b_i\leq i$. For the reason for the terminology
``parking function,'' as well as additional results and references,
see \cite[Exercise~5.49]{stanleyv2}. A basic result of Konheim and
Weiss is that the number of parking functions of length $n$ is
$(n+1)^{n-1}$.  

Write $\pa(\dd)$ for the set of all parking functions of length $\dd$.
For $\xx=(x_1,\dots,x_\dd)\in\mathbb{N}^\dd$ define an
$\xx$-\emph{parking function} to be a sequence $(a_1,\dots,a_\dd)$ of
positive integers whose decreasing rearrangement $b_1\leq\cdots\leq
b_\dd$ satisfies $b_i\leq x_1 + \cdots + x_i$. Thus an ordinary
parking function corresponds to the case $\xx=(1,1,\dots, 1)$. Let
$P_\dd(\xx)$ denote the number of $\xx$-parking functions. Note that
$P_\dd(\xx)=0$ if $x_1=0$.

\begin{thm}{pux}
  \bea P_\dd(\xx) = 
  \sum_{(a_1,\dots,a_\dd)\in\pa(\dd)} x_{a_1}\cdots x_{a_\dd} 
  =
  n! V_\dd(\xx)
  \label{pux1}
    \eea
\end{thm}

\proof
Given $(a_1,\dots,a_\dd)\in\pa(\dd)$, replace each $i$ by an integer
in the set $\{ x_1+\cdots+x_{i-1}+1,\dots,x_1+\cdots+x_i\}$. The
number of ways to do this is given by the middle expression in
(\ref{pux1}), and every $\xx$-parking function is obtained exactly
once in this way.  This yields the first equality. The second equality
follows from the expansion \re{vform} of $V_\dd(\xx)$, since a parking
function is obtained by choosing $\kk\in K_\dd$, forming a sequence
with $k_i$ $i$'s, and permuting its elements in ${\dd\choose
k_1,\dots,k_\dd}$ ways.
\endpf

Take $x_i = 1$ for all $i$ in \re{pux1} and apply 
\re{abform} for $a = b = 1$ to recover the result of
\cite{kon-weiss} that the number of parking
functions of length $\dd$ is $(\dd+1)^{\dd-1}$. 
We note that formula \re{abform} can be given a simple combinatorial
proof generalizing the proof of Pollak \cite[p.\ 13]{fo-ri} for the
case of ordinary parking functions; see \cite[p.\ 10]{rs:nc} for the
case $a=b$.  We note that Theorem \ref{pux} also gives enumerative
interpretations of formulae \re{abc} and \re{abcd1}. Presumably these
formulae too could be derived combinatorially in the setting of
parking functions, but we will not attempt that here.

An interesting special case of Theorem~\ref{pux} arises when we take
$x_i=q^{i-1}$ for some $q>0$. In this case we have
  $$ \dd!\,
     V_\dd(1,q,q^2,\dots,q^{\dd-1})=\sum_{(a_1,\dots,a_\dd)\in
     \pa(\dd)} q^{a_1+\cdots +a_\dd-\dd}. $$
It follows from a result of Kreweras \cite{kreweras2} (see also
\cite[Exer.\ 5.49(c)]{stanleyv2}) that also
  $$ \dd!\, V_\dd(1,q,q^2,\dots,q^{\dd-1}) = q^{{\dd\choose
      2}}I_\dd(1/q), $$ 
where $I_\dd(q)$ is the \emph{inversion enumerator of labeled trees}.

We can generalize equation (\ref{abform}) by giving a simple product
formula for the Ehrhart polynomial $i(\Pi_\dd(\xx),r)$ of
$\Pi_\dd(\xx)$ in the case $\xx= (a,b,b,\dots,b)$ (see
Theorem~\ref{latprod}).  First we need to discuss another way to
interpret $N(\Pi_\dd(\xx))$.

Let $\lambda=(\lambda_1,\dots,\lambda_\ell)$ be a partition, so
$\lambda_i\in \mathbb{N}$ and $\lambda_1\geq \cdots\geq
\lambda_\ell\geq 0$. A \emph{plane partition} of \emph{shape}
$\lambda$ and \emph{largest part at most} $m$ is an array $\pi=
\left(\pi_{ij}\right)$ of integers $1\leq \pi_{ij}\leq m$, defined for
$1\leq i\leq\ell$ and $1\leq j\leq \lambda_i$, which is weakly
decreasing in rows and columns. For instance, the plane partitions of
shape $(2,1)$ and largest part at most $2$ are given by
  $$ \begin{array}{lclclclcl} 11 & & 21 & & 22 & & 21 & & 22\\[-.3em]
                              1  & & 1  & & 1  & & 2  & & 2 
   \end{array}, $$
where we only display the positive parts $\pi_{ij}>0$.
Basic information on plane partitions may be found
in \cite[{\S}{\S}7.20--7.22]{stanleyv2}. If  $\xx=(x_1,\dots,x_\dd)$
$\,\in \mathbb{N}^\dd$ then set 
$$\uu=(u_1,\dots,u_\dd)=
(x_1,x_1+x_2,\cdots,x_1+\cdots+x_\dd)$$ and write
$\tilde{\uu}=(u_\dd, \dots,u_1)$, so that $\tilde{\uu}$ is a partition. 

\begin{thm} {pp}
Let $\xx\in\mathbb{N}^\dd$. Then $N(\Pi_\dd(\xx))$ is equal to the number
of plane partitions of shape $\tilde{\uu}$ and largest part at most $2$.
\end{thm}

\proof
If $(y_1,\dots,y_\dd)\in\Pi_\dd(\xx)\cap \mathbb{Z}^\dd$, then define
the plane partition $\pi$ of shape $\uu$ to have $y_1+\cdots+y_i$ 2's
in row $\dd+1-i$ and the remaining entries equal to 1. This sets up a
bijection between the integer points in $\Pi_\dd(\xx)$ and the plane
partitions of shape $\tilde{\uu}$ and largest part at most 2.
\endpf

\textsc{Note.} Because of the connection given by Theorem~\ref{pp}
between integer points in $\Pi_\dd(\xx)$ and plane partitions, a
number of results concerning $\Pi_\dd(\xx)$ appear already (sometimes
implicitly) in the plane partition literature. In particular, consider
the determinantal formula (\ref{steck2}) of Steck. 
Let $j'_i=j_i-i$,
$b'_i=b_i-i+1$, and $c'_i=c_i-i-1$. We are then counting sequences
$j'_1\leq j'_2\leq \cdots\leq j'_\dd$ satisfying $b'_i\leq j'_i\leq
c'_i$. If $b'_i>b'_{i+1}$ then we can replace $b'_{i+1}$ by $b'_i$
without affecting the sequences $j'_1\leq \cdots\leq j'_\dd$ being
counted. Similarly if $c'_i>c'_{i+1}$ we can replace $c'_i$ with
$c'_{i+1}$. Moreover, clearly the number of sequences being counted is
not changed by adding a fixed integer $k$ to each $b'_i$ and $c'_i$.
Hence it costs nothing to assume that $0\leq b'_1\leq\cdots\leq
b'_\dd$ and $0\leq c'_1\leq \cdots\leq c'_\dd$ (with $b'_i\leq c'_i$).
Let $\lambda=(c'_\dd,\dots,c'_1)$ and $\mu =(b'_\dd, \dots,b'_1)$.
Then $\lambda$ and $\mu$ are partitions, and $\mu\subseteq \lambda$ in
the sense of containment of diagrams (see \cite[{\S}7.2]{stanleyv2}).
Let $Y$ denote the poset (actually a distributive lattice) of all
partitions of all nonnegative integers, ordered by diagram
containment. The lattice $Y$ is just \emph{Young's lattice} mentioned
above. In terms of Young's lattice, we see that that the number
$\#(b,c)$ of equation (\ref{steck2}) is just the number of elements
$(j'_\dd,\dots,j'_1)$ in the interval $[\mu,\lambda]$ of $Y$. 
Alternatively, $\#(b,c)$ is the number of multichains
$\mu=\lambda^0\leq \lambda^1\leq \lambda^2=\lambda$ of length two in
the interval $[\mu,\lambda]$ of $Y$. Kreweras
\cite[{\S}2.3.7]{kreweras} gives a determinantal formula for the
number of multichains of any fixed length $k$ in the interval
$[\mu,\lambda]$. (See also 
\cite[Exer.\ 3.63]{stanley86}.) Such a multichain is easily seen to be
equivalent to a plane partition of shape $\lambda/\mu$ with largest
part at most $k$. When specialized to $k=2$, Kreweras' formula becomes
precisely our equation~(\ref{steck3}). Moreover, the special case
$\mu=\emptyset$ of Kreweras' formula was already known to MacMahon
(put $x=1$ in the implied formula for $GF(p_1p_2\cdots p_m;n)$ in
\cite[p.\ 243]{macmahon}). By Theorem~\ref{pp} the number of elements
of the interval $[\emptyset,\lambda]$ is just
$N(\Pi_\dd(\xx))$, where $\lambda$ is the partition $\tilde{\uu}$ of
Theorem~\ref{pp}. Hence in some sense MacMahon already knew a
determinantal formula for $N(\Pi_\dd(\xx))$ and thus also (by taking
leading coefficients of $N(\Pi_\dd(r\xx))$ regarded as a polynomial in
$r$) for the volume $V_\dd(\xx)$. 

\begin{thm} {latprod}
Let $a,b\in\mathbb{N}$ and $\xx=(a,b,b,\dots,b)\in\mathbb{N}^\dd$. Then
the Ehrhart polynomial $i(\Pi_\dd(\xx))$ is given by
  \beq i(\Pi_\dd(\xx),r) = \frac{1}{\dd!} (ra+1)(r(a+\dd b)+2)
    (r(a+\dd b)+3) \cdots (r(a+\dd b)+\dd). \label{ehrhartpi} \eeq
In particular, the number $N(\Pi_\dd(\xx))$ of integer points in
$\Pi_\dd(\xx)$ satisfies
  $$ N(\Pi_\dd(\xx)) = \frac{1}{\dd!} (a+1)(a+\dd b+2)(a+\dd b+3) \cdots
     (a+\dd b+\dd). $$
\end{thm}

\textbf{First proof.} The theorem is simply a restatement of a
standard result in the subject of ballot problems and lattice path
enumeration, going back at least to Lyness \cite{lyness}, and with
many proofs. A good discussion appears in
\cite[{\S}{\S}1.4--1.6]{mohanty}. See also \cite[{\S}1.3,\
Lemma~3B]{narayana}. 

\textbf{Second proof.}
We give a proof different from the proofs alluded to above, because it
has the virtue of generalizing to give Theorem~\ref{latprod2} below.
The polytope $r\Pi_\dd(\xx)$ is just $\Pi_\dd(r\xx)$. Hence by
Theorem~\ref{pp} $i(\Pi_\dd(\xx),r)$ is just the number of plane
partitions of shape $r\uu$ and largest part at most 2. Identify the
partition $\uu$ with its \emph{diagram}, consisting of all pairs
$(i,j)$ with $1\leq i\leq \dd$ and $1\leq
j\leq\tilde{u}_i=a+(\dd-i)b$. Define the \emph{content} $c(s)$ of
$s=(i,j)\in\tilde{\uu}$ by $c(s)=j-i$ (see \cite[p.\
373]{stanleyv2}). An explicit formula for the 
number of plane partitions of shape $\uu$ and \emph{any} bound on the
largest part was first obtained by Proctor and is
discussed in \cite[Exer.\ 7.101]{stanleyv2} (as well as a
generalization due to Krattenthaler). Proctor's formula for the case
at hand gives
  $$ i(\Pi_\dd(\xx),r) = \prod_{{s=(i,j)\in r\tilde{\tuu}\atop
      \dd+c(s)\leq 
      r\tilde{u}_i}} \frac{1+\dd+c(s)}{\dd+c(s)}\
    \prod_{{s=(i,j)\in r\tilde{\tuu}\atop \dd+c(s)>
      r\tilde{u}_i}} \frac{rb+1+\dd+c(s)}{\dd+c(s)}. $$
When all the factors of the above products are written out, there is
considerable cancellation. The only denominator factors that survive are
those indexed by $(i,1)$, $1\leq i\leq\dd$, yielding the denominator
$\dd!$. The surviving numerator factors are $ra+1$ (indexed by
$(\dd,ra)$) and $r(a+\dd b)+k$, $2\leq k\leq \dd$ (indexed by
$(1,r(a+(\dd-1)b)-n+k)$), the last $\dd-1$ squares in the first row of
$\tilde{\uu}$).   
\endpf

Note from (\ref{ehrhartpi}) that the leading coefficient of
$i(\Pi_\dd(\xx),r)$ (and hence the volume $V_\dd(\xx)$ of $\Pi_\dd(\xx)$) is
given by $a(a+\dd b)^{\dd-1}$, agreeing with equation (\ref{abform}).

There is a straightforward generalization of Theorems~\ref{pp} and
\ref{latprod} involving plane partitions of shape $\uu$ with largest
part at most $m+1$ (instead of just $m+1=2$). Given $\xx\in\mathbb{N}^\dd$
as before, let $\Pi_\dd^m(\xx)\subset\reals^{\dd m}$ be the polytope of
all $\dd\times m$ matrices $\left( y_{ij}\right)$ satisfying
$y_{ij}\geq 0$ and
  $$ v_{i1}\leq v_{i2}\leq\cdots\leq v_{im}\leq x_1+\cdots+x_i, $$
for $1\leq i\leq \dd$, where
  $$ v_{ij} = y_{i1}+y_{i2}+\cdots+y_{ij}. $$
Thus $\Pi_\dd^1(\xx)=\Pi_\dd(\xx)$. Then the proof of Theorem~\ref{pp}
carries over \emph{mutatis mutandis} to show that $N(\Pi_\dd^m(\xx))$ is
the number of plane partitions of shape $\tilde{\uu}$ and largest part
at most $m+1$. The result of Proctor mentioned above gives an explicit
formula for this number when $\xx=(a,b,b,\dots,b)$. Replacing $\xx$ by
$r\xx$ and computing the leading coefficient of the resulting polynomial
in $r$ gives a formula for the volume $V_\dd^m(\xx)$ of
$\Pi_\dd^m(\xx)$. This computation is similar to that in the proof of
Theorem~\ref{latprod}, though the details are more complicated. We
merely state the result here without proof. Is there a direct
combinatorial proof similar to the proofs of Theorem~\ref{latprod}
(the case $m=1$ of Theorem~\ref{latprod2}) appearing in \cite{mohanty}
and \cite{narayana}? 

\begin{thm} {latprod2}
Let $\xx=(a,b,b,\dots,b)\in\mathbb{N}^\dd$. Then
  $$ (\dd m)!\,V_\dd^m(\xx) = 1!\,2!\,\cdots m!\,f^{\langle m^\dd\rangle}
      (\dd+m)^{\dd-1}(\dd+m-1)^{\dd-2}\cdots(\dd+1)^{\dd-m}, $$
where $f^{\langle m^\dd\rangle}$ denotes the number of standard Young
tableaux of shape $\langle m^\dd\rangle=$ $(m,m,\dots,m)$ ($\dd$ $m$'s
in all), given explicitly by the ``hook-length formula'' \cite[Cor.\
7.21.6]{stanleyv2}. 
\end{thm}

\medskip
There is a further generalization of the polytope $\Pi_\dd(\xx)$ which
deserves mention. Let $\xx=(x_1,\dots,x_\dd)\in\nreals^\dd$ and
$\zz=(z_1,\dots,z_\dd)\in\nreals^\dd$, with $v_i=z_1+\cdots+z_i\leq
x_1+\cdots+x_i=u_i$. Let $\Pi_\dd(\zz,\xx)$ be the polytope of all points
$(y_1,\dots,y_\dd)\in \reals^\dd$ satisfying
  $$ y_i\geq 0,\ \mbox{for}\ 1\leq i\leq n $$
  $$ v_i\leq y_1+\cdots +y_i\leq u_i. $$
Thus $\Pi_\dd(\xx)=\Pi_\dd(\bm{0},\xx)$. Much of the theory of $\Pi_\dd(\xx)$
extends to $\Pi_\dd(\zz,\xx)$. Rather than enter into the details here, we
simply illustrate with the case $\dd=2$ how the polyhedral
decomposition of $\Pi_\dd(\xx)$ with chambers $\Lambda_\pi$ extends to
$\Pi_\dd(\zz,\xx)$. In general, the chambers $\Lambda_\pi$ of a decomposition
of $\Pi_\dd(\zz,\xx)$ into a product of simplices will be obtained from
linear extensions $\pi=a_1 a_2 \cdots a_{3\dd}$ of
$\bm{3}\times\bm{\dd}$. Let the elements of $\bm{3}\times \bm{\dd}$ be
$\alpha_1, \dots,\alpha_{3\dd}$ with $\alpha_1<\cdots<\alpha_\dd$,
$\alpha_{\dd+1}<\cdots<\alpha_{2\dd}$,
$\alpha_{2\dd+1}<\cdots<\alpha_{3\dd}$, and $\alpha_i<\alpha_{\dd+i}
<\alpha_{2\dd+i}$ for $1\leq i\leq n$. Then $\pi$ corresponds to the
chamber
  \beq 0\leq f(\alpha_1)\leq\cdots\leq f(\alpha_{3\dd}), 
   \label{3nfacet} \eeq
where 
  $$ f(\alpha_i) = \left\{ \begin{array}{rl}
     v_i, & \mbox{if}\ 1\leq i\leq \dd\\
     y_1+\cdots+y_{i-\dd}, & \mbox{if}\ \dd+1\leq i\leq 2\dd\\
     u_{i-2\dd}, & \mbox{if}\ 2\dd+1\leq i\leq 3\dd.
   \end{array} \right. $$
There is one important difference between this decomposition and the
analogous one for $\Pi_\dd(\xx)$, namely, in the present case some of
the chambers $\Lambda_\pi$ will actually be \emph{empty} and should be
ignored. (Of course $\Lambda_\pi$ isn't really a chamber if it's empty.)
The question of which are empty will depend on the relative order of
the numbers $u_1,\dots,u_\dd$ and $v_1,\dots,v_\dd$. In the ``generic''
case when each $x_i>0$ and $z_i>0$ there are $C_\dd$ (a Catalan number)
relative orderings of the $u_i$'s and $v_i$'s (since $u_1<\cdots<u_\dd$,
$v_1<\cdots<v_\dd$, and $v_i\leq u_i$). More generally, we can change
some of the $\leq$ signs in equation (\ref{3nfacet}) to $<$ signs, in
accordance with the descents of the corresponding linear extension
$\pi$, so that we obtain a decomposition of $\Pi_\dd(\zz,\xx)$ into
pairwise disjoint cells from which we can compute the lattice point
enumerator $N(\Pi_\dd(\zz,\xx))$.

Let us illustrate the above discussion in the case $\dd=2$. The linear
extensions of $\bm{3}\times\bm{2}$, using the labeling just described,
are given by
  $$ \begin{array}{l} 1\,2\,3\,4\,5\,6\\ 1\,2\,3\,5\,4\,6\\
       1\,3\,2\,4\,5\,6\\ 1\,3\,2\,5\,4\,6\\ 1\,3\,5\,2\,4\,6
     \end{array}. $$ 
Thus the following sets (possibly empty) give a decomposition of
$\Pi_\dd(\zz,\xx)$ into pairwise disjoint cells:
  $$ \begin{array}{l} 0\leq v_1\leq v_2\leq y_1\leq y_1+y_2
       \leq u_1\leq u_2\\
      0\leq v_1\leq v_2\leq y_1\leq u_1<y_1+y_2\leq u_2\\
      0\leq v_1\leq y_1<v_2\leq y_1+y_2\leq u_1\leq u_2\\
      0\leq v_1\leq y_1<v_2\leq u_1<y_1+y_2\leq u_2\\
      0\leq v_1\leq y_1\leq u_1<v_2\leq y_1+y_2\leq u_2.
   \end{array} $$
The first four cells are nonempty provided $v_2\leq u_1$, while the
last cell is nonempty provided $v_2>u_1$. Hence we read off that
  $$  N(\Pi_\dd(\zz,\xx)) = \left\{ \begin{array}{rl} A, & x_1\geq
      z_1+z_2\\ \tbbc{x_1-z_1+1}{1}\tbbc{x_1+x_2-z_1-z_2+1}{1}, & 
      x_1<z_1+z_2, \end{array} \right. $$
where
  $$ A= \bbc{x_1-z_1-z_2+1}{2} +\bbc{x_1-z_1-z_2+1}{1}\bbc{x_2}{1} $$
\vspace{-.3em}
   $$ \qquad + \bbc{z_2}{1}\bbc{x_1-z_1-z_2+1}{1}
   +\bbc{z_2}{1}\bbc{x_2}{1}. $$  

\section{A subdivision of $\Pi_\dd(\xx)$ connected with the associahedron}
\label{sec:assoc}
In this section we describe a polyhedral subdivison
$(\hat{\Pi}_\dd(\kk;\xx),\ \kk\in K_\dd)$ of
$\Pi_\dd(\xx)$ different from the subdivision discussed in
Section~\ref{ordercone}. This subdivision is closely related to a
convex polytope known as the \emph{associahedron}, defined as
follows. Let $E_{\dd+2}$ be a convex $(\dd+2)$-gon. A \emph{polygonal
decomposition} of $E_{\dd+2}$ 
consists of a set of diagonals of $E_{\dd+2}$ that do not cross in
their interiors. Hence the maximal polygonal decompositions are the
triangulations, and contain exactly $\dd-1$ diagonals. Let
$\mathrm{dec}(E_{\dd+2})$ denote the poset of all polygonal
decompositions of $E_{\dd+2}$, ordered by inclusion, with a top
element $\hat{1}$ adjoined. It 
was first shown by C. W. Lee \cite{lee89} and M. Haiman \cite{haiman}
that $\mathrm{dec}(E_{\dd+2})$ is the face lattice of an
$(\dd-1)$-dimensional convex polytope ${\cal A}_\dd$, known as the
\emph{associahedron} or \emph{Stasheff polytope}. (Earlier Stasheff
\cite{stasheff} defined the dual of the associahedron as a simplicial
complex and constructed a geometric realization as a convex body but
not as a polytope.)  A vast generalization is discussed in \cite[Ch.\ 
7]{g-k-z}. For some further information see \cite[Exer.\ 
6.33]{stanleyv2}.

We next give a somewhat different description of the associahedron (or
more precisely, of its face lattice) that is most convenient for our
purposes. A \emph{fan} in $\reals^m$ is a (finite) collection $\bm{F}$
of pointed polyhedral cones (with vertices at the origin) satisfying
the two conditions:
  \begin{itemize} \item If ${\cal C},{\cal C}'\in\bm{F}$ then ${\cal
C}\cap {\cal C}'$ is a face (possibly consisting of just the origin) of
${\cal C}$ and ${\cal C}'$.
  \item If ${\cal C}\in\bm{F}$ and ${\cal C}'$ is a face of ${\cal
      C}$, then ${\cal C}'\in\bm{F}$.
  \end{itemize}
A fan $\bm{F}$ is called \emph{complete} if $\bigcup_{{\cal
C}\in\sbm{F}} = \reals^m$.

We will define a fan whose chambers are indexed by plane binary trees
with $n$ internal vertices. The definition of a plane tree may be
found for instance in \cite[Appendix]{stanley86}. The key point is
that the subtrees of any vertex are linearly ordered $T_1,\dots,T_k$,
indicated in drawing the tree (with the root on the bottom) by placing
the subtrees in the order $T_1,\dots,T_k$ from left to right. A
\emph{binary} plane tree is a plane tree for which each vertex $v$ has
zero or two subtrees. In the latter case we call the vertex an
\emph{internal} vertex. Otherwise $v$ is a \emph{leaf} or
\emph{endpoint}. We will always regard plane trees as being drawn with
the root at the bottom. 

Let $T$ be a plane binary tree with $\dd$ internal vertices (so
$\dd+1$ leaves). The number of such trees is the Catalan number $C_n$
\cite[6.19(d)]{stanleyv2}. 
Do a depth-first search through $T$ (as defined e.g.\ 
in \cite[pp.\ 33--34]{stanleyv2}) and label the internal vertices
$1,2,\dots,\dd$ in the order they are first encountered \emph{from
above}. Equivalently, every internal vertex is greater than those in
its left subtree, and smaller than those in its right subtree.  We
call this labeling of the internal vertices of $T$ the \emph{binary
search labeling}.  Figure~\ref{fig:ptree} gives an example when
$\dd=4$. Let $y_1,\dots,y_{\dd-1}$ denote the coordinates in
$\reals^{\dd-1}$. If the internal vertex $i$ of $T$ (using the
labeling just defined) is covered by $j$ and $i<j$, then associate
with the pair $(i,j)$ the inequality
  \beq y_{i+1}+y_{i+2}+\cdots+y_j\leq 0, \label{edgefacet1} \eeq
while if $i>j$ then associate with $(i,j)$ the inequality
  \beq y_{j+1}+y_{j+2}+\cdots+y_i\geq 0. \label{edgefacet2} \eeq
We get a system of $\dd-1$ homogeneous linear inequalities that define
a simplicial
cone ${\cal C}_T$ in $\reals^{\dd-1}$. For example, the inequalities
corresponding to the tree of Figure~\ref{fig:ptree} are given by
  \beas y_2 & \leq & 0\\ y_2+y_3 & \geq & 0\\ y_4 & \geq & 0.
  \eeas

 \setlength{\unitlength}{.6pt} 
\begin{figure}
   \begin{picture}(200,170)(-200,-10)
  \put(100,0){\circle*{5}}
  \put(50,50){\circle*{5}}
  \put(150,50){\circle*{5}}
  \put(10,100){\circle*{5}}
  \put(90,100){\circle*{5}}
  \put(110,100){\circle*{5}}
  \put(190,100){\circle*{5}}
  \put(50,150){\circle*{5}}
  \put(130,150){\circle*{5}}
  \put(100,0){\line(-1,1){50}}
  \put(100,0){\line(1,1){50}}
  \put(50,50){\line(-4,5){40}}
  \put(50,50){\line(4,5){80}}
  \put(90,100){\line(-4,5){40}}
  \put(150,50){\line(-4,5){40}}
  \put(150,50){\line(4,5){40}}
  \put(80,-7){3}
  \put(32,42){1}
  \put(72,95){2}
  \put(158,42){4}
  \end{picture}
\caption{A plane tree with the binary search labeling of its internal
  vertices} 
\label{fig:ptree}
\end{figure}
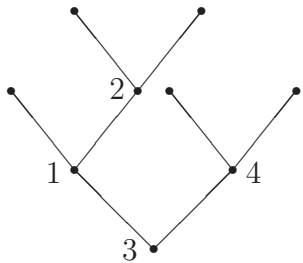

\begin{figure}
\centerline{\epsfig{figure=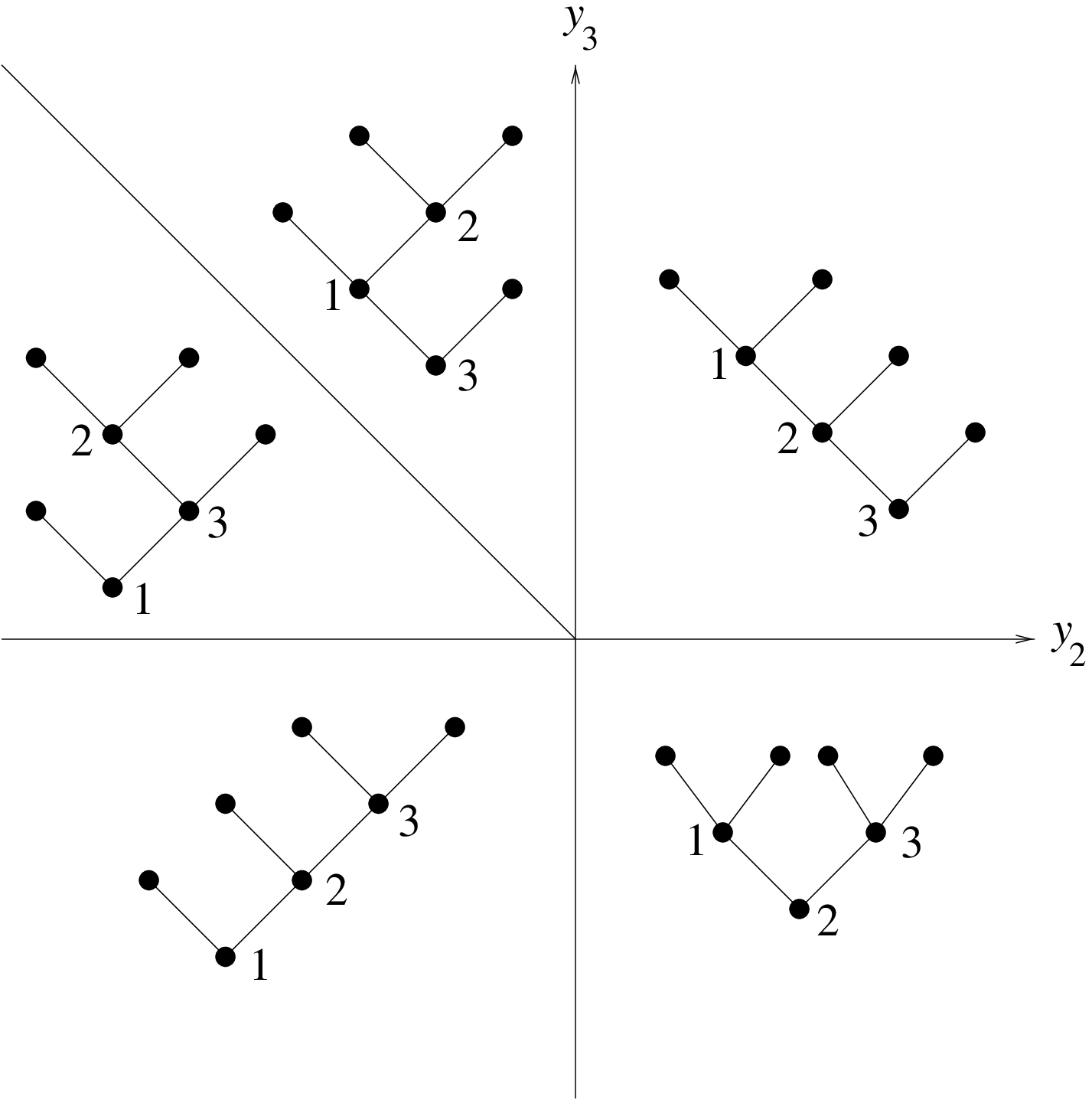}}
\caption{The fan $\bm{F}\!_3$}
\label{fig:fan}
\end{figure}

It is not hard to check that these $C_\dd$ cones, as $T$ ranges over all
plane binary trees with $\dd$ internal vertices, form the chambers of a
complete fan $\bm{F}\!_\dd$ in $\reals^{\dd-1}$.  For instance,
Figure~\ref{fig:fan} shows the fan $\bm{F}\!_3$.

\begin{thm}{associahedron}
The face poset $P(\bm{F}\!_\dd)$ of the fan $\bm{F}\!_\dd$, with a top
element $\hat{1}$ adjoined, is isomorphic to the dual
$\mathrm{dec}(E_{\dd+2})^*$ of the face lattice of the associahedron
${\cal A}_{\dd+2}$. 
\end{thm}

\begin{figure}
\centerline{\psfig{figure=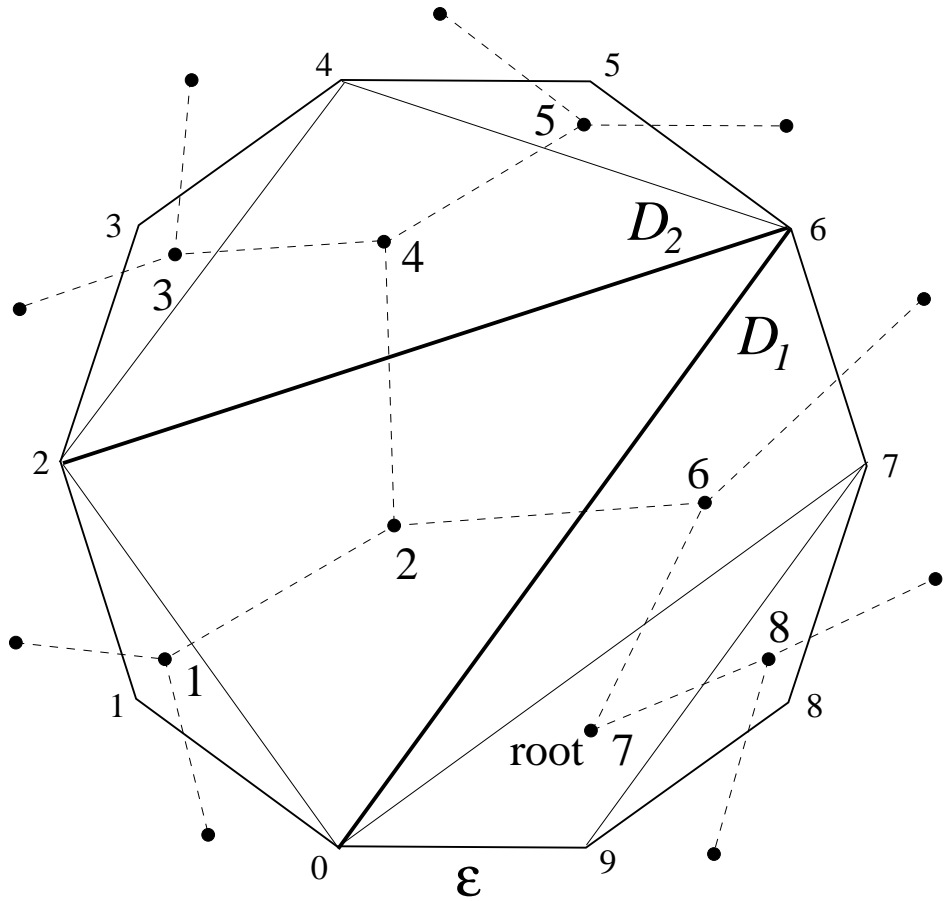}}
\caption{A triangulated 10-gon and the corresponding plane binary tree
  $T$}
\label{triang}
\end{figure}

\proof The face lattice of a complete fan is completely determined by
the incidences between the chambers and extreme rays. (See \cite[Exer.\ 
3.12]{stanley86} for a stronger statement.) The chambers of $\bm{F}\!_\dd$
have already been described in terms of plane binary trees. There is a
well-known bijection between plane binary trees on $2\dd+1$ vertices and
triangulations of a convex $(\dd+2)$-gon $E_{\dd+2}$. This bijection is
explained for instance in \cite[Cor.\ 6.2.3]{stanleyv2}. In
particular, to define the bijection we first need to fix an edge
$\varepsilon$ of $E_{\dd+2}$, called the \emph{root edge}. We hope that
Figure~\ref{triang} will make this bijection clear; see
the previous reference for further details. Thus we have a
bijection between the chambers ${\cal C}$ of $\bm{F}\!_\dd$ and the
triangulations of the convex $(\dd+2)$-gon $E_{\dd+2}$.

We now describe the extreme rays $R$ of $\bm{F}\!_\dd$. We can
describe $R$ uniquely by specifying one nonzero point on $R$. We will
index these points by the diagonals $D$ of a convex $(\dd+2)$-gon
$E_{\dd+2}$. Label the vertices of $E_{\dd+2}$ as $0,1,\dots,\dd+1$
clockwise beginning
with one vertex of $\varepsilon$ and ending with the other. Let $e_i$
denote the unit coordinate vector corresponding to the coordinate
$y_i$ in the space $\reals^{\dd-1}$ with coordinates
$y_2,\dots,y_\dd$. Given the diagonal $D$ between vertices $i<j$ of
$E_{\dd+2}$, associate a point $p_D\in \reals^{\dd-1}$ as follows:
  $$ p_D= \left\{ \begin{array}{rl} e_j, & \mbox{if}\ i=0\\
    -e_{i+1}, & \mbox{if}\ j=\dd+1\\
    e_j-e_{i+1}, & \mbox{otherwise}. \end{array} \right. $$
We claim that the ray $\{\alpha p_D\,:\,\alpha\in\nreals\}$ is the
extreme ray of $\bm{F}\!_\dd$ that is the intersection of all the chambers
of $\bm{F}\!_\dd$ corresponding to the triangulations of $E_{\dd+2}$
that contain $D$. From this claim the proof of the theorem follows
(using the fact that $\bm{F}\!_\dd$ is a \emph{simplicial} fan, i.e.,
every face is a simplicial cone).

Consider first the diagonal $D$ with vertices 0 and $j$. Let
$\Upsilon$ be a triangulation of $E_{\dd+2}$ containing $D$. The
internal vertices of $T$ corresponding to the regions (triangles) of
the triangulation $\Upsilon$. Because of our procedure for labeling
the internal vertices of a plane binary tree $T$, it follows that the
labels of the internal vertices ``above'' $D$ (i.e., on the opposite
side of $D$ as the root edge $\varepsilon$) will be $1,2,\dots,j-1$,
while the internal vertices below $D$ will be labeled $j,j+1,\dots,
\dd$. (See Figure~\ref{triang} for an example with $\dd=8$. The
diagonal $D$ in question is labeled $D_1$ and connects vertex 0 to
vertex $j=6$. The plane binary tree $T$ is drawn with dashed lines.)
Consider the internal edges of $T$ that give rise (\emph{via}
equations (\ref{edgefacet1}) and (\ref{edgefacet2})) to chambers whose
equations involve $y_j$. No such edge can appear below $D$, since $j$
is the least vertex label appearing below $D$. Similarly no such edge
can appear above $D$, since only vertices less than $j$ appear above
$D$. Hence such an edge must cross $D$. The top (farthest from the
root) vertex $a$ of this edge is $<j$, while the bottom vertex $b$ is
$\geq j$. Hence the chamber equation is given by
$y_{a+1}+y_{a+2}+\cdots+y_b\geq 0$, where $a< j$ and $b\geq j$. Hence
the point $e_j$ lies on this chamber, and so the ray through $e_j$ is
the intersection of the chambers corresponding to triangulations
containing $D$.

A completely analogous argument holds for the diagonal $D$ with
vertices $i$ and $\dd+1$.

Finally suppose that $D$ has vertices $i,j$ where $0<i<j<\dd+1$. The
internal vertices of $T$ appearing above $D$ will be labeled
$i+1,i+2,\dots, j-1$, while the remaining vertex labels appear below
$D$. (See Figure~\ref{triang}, where the diagonal $D$ in question is
labeled $D_2$, and where $i=2$ and $j=6$.) Consider an internal edge
of $T$ whose vertex labels are $a$ and $b$ where $a\leq i$ and
$i+1\leq b<j$. These are precisely the edges whose corresponding chamber
equation (either $y_{a+1}+y_{a+2}+\cdots+y_b\geq 0$ or
$y_{a+1}+y_{a+2}+\cdots+y_b\leq 0$) involves $y_{i+1}$ but not $y_j$.
Since $b$ appears above $D$ and $a$ below, the chamber equation is in
fact $y_{a+1}+y_{a+2}+\cdots+y_b\leq 0$. In particular, the point
$e_j-e_{i+1}$ lies on the chamber. Similarly, consider an internal edge
of $T$ whose labels are $a$ and $b$ where $i+1\leq a<j$ and $j\leq b$.
These are precisely the edges whose corresponding chamber equation
(again either $y_{a+1}+y_{a+2}+\cdots+y_b\geq 0$ or
$y_{a+1}+y_{a+2}+\cdots+y_b\leq 0$) involves $y_j$ but not $y_{i+1}$.
Since $b$ appears below $D$ and $a$ above, the chamber equation is in
fact $y_{a+1}+y_{a+2}+\cdots+y_b\geq 0$. In particular, the point
$e_j-e_{i+1}$ lies on the chamber. Every other chamber equation either
involves neither $y_{i+1}$ nor $y_j$, or else involves both (with a
coefficient 1). Hence $e_{i+1}-e_j$ lies on every chamber corresponding
to a triangulation containing $D$, so the intersection of these chambers
is the ray containing $e_j-e_{i+1}$. This completes the proof of the
claim, and with it the theorem.
\endpf

The connection between $\Pi_\dd(\xx)$ and the fan $\bm{F}\!_\dd$ is
provided by the concept of a plane tree with edge lengths.  If we
associate with each edge $e$ of the plane tree $T$ a positive real
number $\ell(e)$, then we call the pair $(T,\ell)$ a \emph{plane tree
with edge lengths}.  Such a tree can be drawn by letting the length
of each edge $e$ be $\ell(e)$.

\begin{figure}
\centerline{\psfig{figure=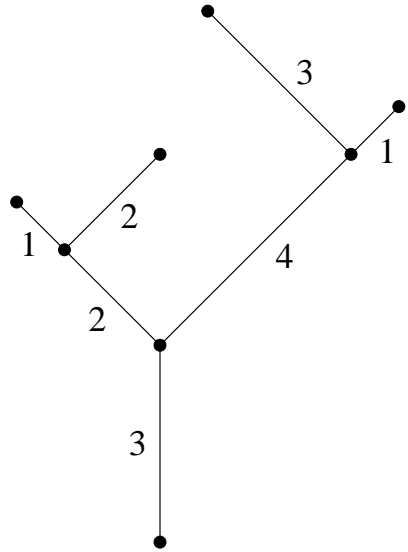}}
\caption{A planted plane binary tree with edge lengths}
\label{eltree}
\end{figure}

Now fix a real number $s>0$, which will be the sum of the edge lengths
of a plane tree. Let $\xx=(x_1,\dots,x_\dd)\in\preals^\dd$ with $\sum
x_i<s$. Let $\yy=(y_1,\dots,y_\dd)\in \preals^\dd$ with
$y_1+\cdots+y_i\leq x_1+\cdots +x_i$ for $1\leq i\leq \dd$. We
associate with the pair $(\xx,\yy)$ a plane tree with edge lengths
$\varphi(\xx,\yy)=(\bar{T},\ell)$ as follows. Start at the root and
traverse the tree in preorder (or depth-first order) \cite[pp.\ 
33--34]{stanleyv2}.  First go up a distance $x_1$, then down a
distance $y_1$, then up a distance $x_2$, then down a distance $y_2$,
etc. After going down a distance $y_\dd$, complete the tree by going
up a distance $x_{\dd+1}=s-x_1-\cdots-x_\dd$ and then down a distance
$y_{\dd+1}=s-y_1-\cdots-y_\dd$. Generically we obtain a \emph{planted
  plane binary tree with edge lengths}, i.e, the root has degree one
(or one child), and all other internal vertices have degree two.
Figure~\ref{eltree} shows the planted plane binary tree with edge
lengths associated with $s=16$ and $\xx=(6,2,7)$, $\yy= (1,4,3)$. If
$\bar{T}$ is a planted plane tree, then we let $T$ denote the tree
obtained by ``unplanting'' (uprooting?) $\bar{T}$, i.e., remove from
$\bar{T}$ the root and its unique incident edge $e$ (letting the other
vertex of $e$ become the root of $T$).

\begin{figure}
\centerline{\psfig{figure=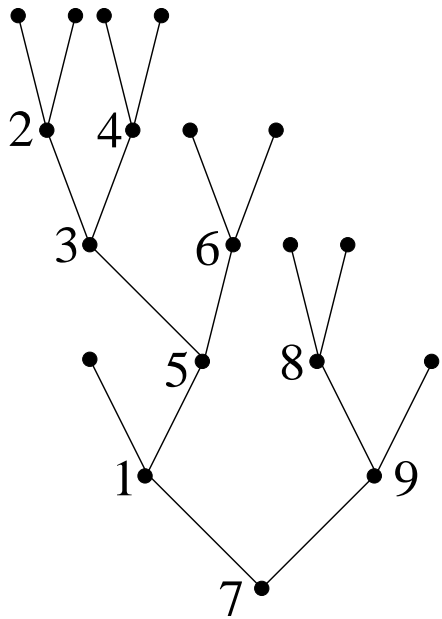}}
\caption{A plane binary tree $T$ with $\kk(T)=(2,3,0,1,0,1,0,2,0)$}
\label{fig:kt}
\end{figure}

Fix the sequence $\xx=(x_1,\dots,x_\dd)$ with $\sum x_i<s$.  For a
plane binary tree $T$ (without edge lengths) with $\dd$ internal
vertices (and hence $\dd+1$ leaves), define $\Delta_T=\Delta_T(\xx)$
to be the set of all $\yy=(y_1,\dots,y_\dd)\in \preals^\dd$ such that
$\varphi(\xx,\yy)=(\bar{T},\ell)$ for some $\ell$. Let ${\cal T}_\dd$
denote the set of plane binary trees with $\dd$ internal vertices. Let
$T\in{\cal T}_\dd$ with the binary search labeling of its internal
vertices as defined earlier in this section. We now define a sequence
$\kk(T) = (k_1,\dots,k_\dd)\in\mathbb{N}^\dd$ as follows: (1) $k_i=0$ if
the left child of vertex $i$ is an internal vertex. (2) If the left
child of vertex $i$ is an endpoint, then let $k_i$ be the largest
integer $r$ for which there is a chain $i=j_1<j_2<\cdots<j_r$ of
internal vertices such that $j_h$ is a left child of $j_{h+1}$ for
$1\leq h\leq r-1$. For instance, if $T$ is the tree of
Figure~\ref{fig:kt} then $\kk(T)=(2,3,0,1,0,1,0,2,0)$.

\begin{lemma} \label{ktbijection}
The map $T\mapsto \kk(T)$ is a bijection from ${\cal T}_\dd$ to the set
$K_\dd$ defined by equation (\ref{kdef}).
\end{lemma}

\proof
Let $\kk(T)=(k_1,\dots,k_\dd)$. The chains $i=j_1<j_2<\cdots<j_r$
described above partition the internal vertices of $T$, so $\sum
k_i=\dd$. Since $k_{j_2}=\cdots=k_{j_r}=0$, it follows that
$k_{h+1}+k_{h+2}+\cdots +k_\dd\leq \dd-h$ for $0\leq h\leq
\dd-1$. Hence $k_1+\cdots+k_h\geq h$, so $\kk(T)\in K_\dd$.

It remains to show that given $\kk=(k_1,\dots,k_\dd)\in K_\dd$, there
is a unique $T\in {\cal T}_\dd$ such that $\kk(T)=\kk$. We can
construct the subtree of internal vertices of $T$ as follows. Let
$T_1$ be defined by starting at the root and making $k_1-1$ steps to
the left. (Each step is from a vertex to an adjacent vertex.) Hence we
have $k_1$ vertices in all, and we are located at the vertex furthest
from the root. Suppose that $T_i$ has been constructed for $i<\dd$,
and that we are located at vertex $v_i$. If $k_{i+1}>0$, then move one
step to the right and $k_{i+1}-1$ steps to the left, yielding the tree
$T_{i+1}$ and the vertex $v_{i+1}$ at which we are located. If
$k_{i+1}=0$, then move down the tree (toward the root) until we have
traversed exactly one edge in a southeast direction. This gives the
tree $T_{i+1}=T_i$ and a new present location $v_{i+1}$. Let
$T=T_\dd$. It is easily checked that the definition of $K_\dd$ ensures
that $T$ is defined (and, though not really needed here, that $v_\dd$
is the root vertex) and $\kk(T)=\kk$. Since there are
$C_\dd=\frac{1}{\dd+1}{2\dd\choose \dd}$ plane binary trees with $\dd$
internal vertices and since $\#K_\dd=C_\dd$, it follows that the map
$T\mapsto\kk(T)$ is a bijection as claimed. (It is also easy to see
directly that $T$ is unique, i.e., if $\kk(T)= \kk(T')$ then $T=T'$.)
\endpf

Now given $t\in\preals$, let $\sigma_k(t)$ denote the $k$-dimensional
simplex of points $(t_1,\dots,t_k)$ satisfying $0\leq t_1\leq t_2 \leq
\cdots\leq t_k\leq t$. Thus
   $$ \mathrm{Vol}(\sigma_k(t)) =\frac{t^k}{k!}. $$ 
By convention
$\sigma_0(t)$ is just a point, with Vol$(\sigma_0(t))=1$. We can now
state the main result of this section.

\begin{thm} {mainassoc}
\textrm{(a)} The sets $\Delta_T(\xx)$, for $T\in{\cal T}_\dd$, form
the maximal faces (chambers) of a polyhedral decomposition $\Gamma_\dd$
of $\Pi_\dd(\xx)$.

\textrm{(b)} Let $\kk(T)=(k_1,\dots,k_\dd)$, where $T\in {\cal T}_\dd$. Then
$\Delta_T(\xx)$ is integrally equivalent (as defined at the beginning
of Section~\ref{pinsec}) to the product
$\sigma_{k_1}(x_1) \times \cdots\times\sigma_{k_\dd}(x_\dd)$, so in
particular
  $$ \mathrm{Vol}(\Delta_T(\xx)) =\frac{x_1^{k_1}}{k_1!}\cdots
     \frac{x_\dd^{k_\dd}}{k_\dd!}. $$

\textrm{(c)} The interior face complex $\Gamma_\dd^\circ$ of
$\Gamma_\dd$ is combinatorially equivalent to the dual associahedron,
i.e., the set of interior faces of $\Gamma_\dd$, ordered by inclusion,
in isomorphic to the face lattice of the dual associahedron.
\end{thm}

\noindent \textbf{Proof of (a).} 
The construction of the plane tree with edge lengths
$\varphi(\xx,\yy)=(\bar{T},\ell)$ is defined if and only if $\yy \in
\Pi_\dd(\xx)$. Since generically $\varphi(\xx,\yy)$ is a planted plane
binary tree, it follows that the sets $\Delta_T(\xx)$, $T\in {\cal
T}_\dd$, form the chambers of a polyhedral decomposition of
$\Pi_\dd(\xx)$. 
\endpf

\noindent \textbf{Proof of (b).} 
Let $\varphi(\xx,\yy)=(\bar{T},\ell)$ as above. Call a vertex $v$
of $\bar{T}$ a \emph{left leaf} if it is a leaf (endpoint) and is the
left child of its parent. Similarly a \emph{right edge} is an edge that
slants to the right as we move away from the root.
Let $P(v)$ be the path from the left leaf
$v$ toward the root that terminates after the first right edge is
traversed (or terminates at the root if there is no such right
edge). Let $c(v)$ be the label of the (internal) vertex covered by
$v$. Then the length of the path $P(v)$ is just $x_{c(v)}$. If
$c(v)=i$, then exactly $k_i$ of the paths $P(u)$ end at the path
$P(v)$. Suppose that these paths are $P(u_1),\dots,P(u_{k_i})$ where
$u_1<\cdots < u_{k_i}$. Then the paths $P(u_j)$ intersect the path
$P(v)$ in the order $P(u_1),\dots,P(u_{k_i})$ from the bottom
up. Hence for each $i$ with $k_i>0$, we can independently place on 
a path of length $x_i$ the $k_i$ points that form the bottoms of the
paths $P(u_j)$. The placement of these points defines a point in a
simplex integrally equivalent to $\sigma_{k_i}(x_i)$, so
$\Delta_T(\xx)$ is integrally equivalent to $\sigma_{k_1}(x_1)
\times\cdots\times \sigma_{k_\dd}(x_\dd)$ as claimed.
\endpf

\begin{figure}
\centerline{\epsfig{figure=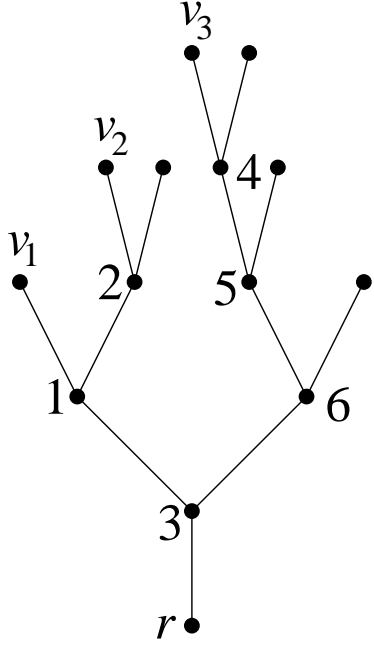}}
\caption{A planted plane binary tree}
\label{fig:treeps}
\end{figure}

\begin{example}
Let $\bar{T}$ be the planted plane binary tree of
Figure~\ref{fig:treeps}. On the path of length $x_1$ from the root $r$
to $v_1$ we can place vertices $1$ and $3$ in bijection with the points of
the simplex $0\leq t_3\leq t_1\leq x_1$ of volume $x_1^2/2$. On the
path of length $x_2$ from 1 to $v_2$ we can place vertex $2$ in
bijection with the points of the simplex $0\leq t_2\leq x_2$, of volume
$x_2$. Finally on the path of length $x_4$ from $3$ to $v_3$ we can
place vertices $4,5,6$ in bijection with the points of the simplex
$0\leq t_6\leq t_5\leq t_4\leq x_4$, of volume $x_4^3/6$. Hence
$\Delta_T$ is integrally equivalent to the product $\sigma_2(x_1)
\times \sigma_1(x_2)\times \sigma_3(x_4)$, of volume
$x_1^2x_2x_4^3/2!\,1!\,3!$. 
\end{example}

It is easy to make the integral equivalence between $\Delta_T$ and
$\sigma_{k_1}(x_1)\times \cdots\times\sigma_{k_\dd}(x_\dd)$ completely
explicit. For instance, in the above example $t_3$ is the distance
between vertices $r$ and 3, so
  $$ t_3 = x_1-y_1+x_2-y_2+x_3-y_3. $$
Similarly,
  $$ t_1 = x_1-y_1. $$
Now $t_2$ is the distance between vertices 1 and 2, so
  $$ t_2 =x_2-y_2. $$
In the same way we obtain
  \beas t_6 & = & x_4-y_y+x_5-y_5+x_6-y_6\\
    t_5 & = & x_4-y_4+x_5-y_5\\
    t_4 & = & x_4-y_4. \eeas

\noindent
\textbf{Proof of (c).} 
Let $\varphi(\xx,\yy)=(\bar{T},\ell)$. Then the
height (or distance from the root) of vertex $i$ is just
$x_1+\cdots+x_i-y_1 -\cdots -y_i=u_i-v_i$. Hence if vertex $i$ is
covered by $j$ then $u_i-v_i< u_j-v_j$. If $i<j$ we get the equation
  \beq (y_{i+1}-x_{i+1})+\cdots+(y_j-x_j)\leq 0, \label{facet1} \eeq
while if $i>j$ we get 
   \beq (y_{j+1}-x_{j+1})+\cdots+(y_i-x_i)\geq 0. \label{facet2} \eeq
Thus these equations, together with $y_i\geq 0$ and
$y_1+\cdots+y_i\leq x_1+\cdots+x_i$, determine $\bar{\Delta}_T$. 

Note that if we replace each $y_k$ by $y_k-x_k$ in the inequalities
(\ref{edgefacet1}) and (\ref{edgefacet2}) defining the chambers of the
fan $\bm{F}\!_\dd$ of Theorem~\ref{associahedron}, then we obtain
precisely the inequalities 
(\ref{facet1}) and (\ref{facet2}). From this we conclude the
following. Given $\xx=(x_1,\dots,x_\dd)\in\nreals^\dd$, translate the
fan $\bm{F}\!_\dd$ so that the center of the translated fan
$\widetilde{\bm{F}}\!_\dd$ is at $(x_2,\dots,x_\dd)$. Add a new $y_1$ axis
and lift $\widetilde{\bm{F}}\!_\dd$ into $\reals^\dd$, giving a
``nonpointed fan'' (i.e., a decomposition of $\reals^n$ satisfying the
definition of a fan except that the cones are nonpointed) which we
denote by $\widetilde{\bm{F}}\!_\dd\times\reals$. (Thus each cone
${\cal 
  C}\in \widetilde{\bm{F}}\!_\dd$ lifts to the nonpointed cone
${\cal C}\times \reals$.) Finally intersect each chamber (maximal cone)
${\cal C}\times\reals$ of $\widetilde{\bm{F}}\!_\dd\times\reals$ with the
polytope $\Pi_\dd(\xx)$. 
Then the polytopes ${\cal C}\cap\Pi_\dd(\xx)$ are just the chambers
$\hat{\Pi}(\kk;\xx)$ of the polyhedral decomposition
${\cal P}_\dd$ of $\Pi_\dd(\xx)$.
Moreover, the interior faces of this decomposition are just the
intersections of \emph{arbitrary} cones in
$\widetilde{\bm{F}}\!_\dd\times\reals$
with $\Pi_\dd(\xx)$. Hence the interior face poset of
${\cal P}_\dd$ is isomorphic to the face
poset of the fan $\bm{F}\!_\dd$, which by Theorem~\ref{associahedron}
is the face lattice of the dual associahedron.
\endpf

\textsc{Notes.} 

The decomposition of $\Pi_\dd(\xx)$ given by
Theorem~\ref{associahedron} is fundamentally different (i.e., has a
different combinatorial type) than that of Theorem~\ref{thm:decomp}.
For instance, when $\dd=3$ Figure~\ref{fig:p3decomp} shows that the
interior face dual complex described by Theorem~\ref{thm:decomp} is
not a decomposition of a convex polytope, unlike the situation in
Theorem~\ref{associahedron}. In that case when $\dd=3$ the interior
face dual complex is just a solid pentagon. The two subdivisions fo
$\Pi_3(\xx)$ are shown explicitly in Figure~\ref{fig3}. 

We are grateful to Victor Reiner for pointing
out to us that Theorem~\ref{associahedron} is related to the
construction of the associahedron appearing 
 in the papers \cite{lee89} and \cite{re-zi}, and that a
$B_\dd$-analogue of this construction
appears in \cite[{\S}3]{bu-re}. Note that the proof of
Theorem~\ref{associahedron} shows that the extreme rays of the fan
$\bm{F}\!_\dd$ are the vectors $e_i$ and $-e_i$ for $1\leq i\leq n-1$,
and $e_i-e_j$ for $1\leq i<j\leq n-1$. As pointed out to us by Reiner,
it follows from \cite{lee89} that we can rescale these vectors (i.e.,
multiply them by suitable positive real numbers) so that their convex
hull is combinatorially equivalent (as defined in the next section) to
the dual associahedron ${\cal A}_{\dd+2}^*$.

Some of the results of this section can be interpreted
probabilistically in terms of the kind of random plane tree with edge
lengths derived from a Brownian excursion by Neveu and Pitman \cite{NP89a}.
It was in fact by consideration of such random trees that we were first
led to the formula \re{vform} for the volume polynomial, with the
geometric interpretation provided by Theorem \ref{mainassoc}.

\section{The face structure of $\Pi_\dd(\xx)$} 

In this section we determine the structure of the faces of
$\Pi_\dd(\xx)$, i.e., a description of the lattice of faces of
$\Pi_\dd(\xx)$ (ordered by inclusion). This description will depend on
the ``degeneracy'' of $\Pi_\dd(\xx)$, i.e., for which $i$ we have
$x_i=0$. Thus let $u_i=x_1+\cdots+x_i$ as usual, and define integers
$1\leq a_1<a_2<\cdots<a_k=\dd$ by
  $$ u_1=\cdots=u_{a_1}<
           u_{a_1+1}=\cdots=u_{a_2}< \cdots<
           u_{a_{k-1}+1}=\cdots=u_{a_k}. $$
We say that two convex polytopes are \emph{combinatorially equivalent}
or have the same \emph{combinatorial type}
if they have isomorphic face lattices. 

\begin{thm} {faces}
Let $a_1,\dots,a_k$ be as above, and set $b_i=a_i-a_{i-1}$ (with
$a_0=0$). Assume (without loss of generality) that $x_1>0$. Then
$\Pi_\dd(\xx)$ is combinatorially equivalent to a product 
$\sigma_{b_1}\times\cdots\times \sigma_{b_k}$, where $\sigma_j$
denotes a $j$-simplex. In particular, if each $x_i>0$ then
$\Pi_\dd(\xx)$ is combinatorially equivalent to an $\dd$-cube.
\end{thm}

\proof
For $1\leq i\leq k$, let ${\cal S}_i=\{ C_{i0},C_{i1},\dots,C_{i,b_i}\}$
denote the set of the following $b_i+1$ conditions $C_{ij}$ on a point
$y\in\Pi_\dd(\xx)$: 
  $$ \begin{array}{rc} \mbox{($C_{i0}$)} & y_{a_{i-1}+1}=y_{a_{i-1}+2} 
     =\cdots=y_{a_i}=0\\[.3em]
    \mbox{($C_{i1}$)} & y_{a_{i-1}+1}=u_i,\ y_{a_{i-1}+2}=y_{a_{i-1}+3}=
        \cdots=y_{a_i}=0\\[.3em]
    \mbox{($C_{i2}$)} & y_{a_{i-1}+2}=u_i,\ y_{a_{i-1}+1}=y_{a_{i-1}+3}=
        \cdots=y_{a_i}=0\\ & \cdots\\[.3em]
    \mbox{($C_{i,b_i}$)} & y_{a_i}=u_i,\ y_{a_{i-1}+1}=y_{a_{i-1}+2}=
        \cdots=y_{a_i-1}=0.
    \end{array} $$
Note that each of the conditions $C_{ij}$ consists of $b_i$
chambers of $\Pi_\dd(\xx)$; we regard $C_{ij}$ as being the set of these
chambers. 
Let $S_i$ denote any subset of ${\cal S}_i$, and let $\cap
S_i=\bigcap_{C\in S_i}C$.  A little thought shows that we can find a
point $y\in\Pi_\dd(\xx)$ lying on all the chambers in each $\cap S_i$, but
not lying on any other chamber of $\Pi_\dd(\xx)$. Moreover, no point of
$\Pi_\dd(\xx)$ can lie on any other collection of chambers of $\Pi_\dd(\xx)$
but on no additional chambers. 

{From} the above discussion it follows that $\Pi_\dd(\xx)$ is
combinatorially equivalent to a product of simplices of dimensions
$b_1,\dots, b_k$, as desired. In particular, $\Pi_\dd(\xx)$ has
$(b_1+1)(b_2+1)\cdots (b_k+1)$ vertices $v$, obtained by choosing $0\leq
j_i\leq b_i$ for each $i$ and defining $v$ to be the intersection of
the chambers in all the $C_{ij_i}$'s.
\endpf

Although $\Pi_\dd(\xx)$ is combinatorial equivalent to a product of
simplices, it is not the case that $\Pi_\dd(\xx)$ is \emph{affinely}
equivalent to such a product. For instance, Figure~\ref{fig2} shows
$\Pi_2(x_1,x_2)$ when $x_1,x_2>0$. We see that $\Pi_2(x_1,x_2)$ is a
quadrilateral and hence combinatorially equivalent to a
square. However, $\Pi_2(x_1,x_2)$ is not a parallelogram and hence not
affinely equivalent to a square. Similarly Figure~\ref{fig3} shows that
$\Pi_3(x_1,x_2,x_3)$ is combinatorially equivalent but not affinely
equivalent to a 3-cube when each $x_i>0$.

\end{document}